\documentclass[12pt,centertags,oneside]{amsart}
\usepackage{amsmath,amstext,amsthm,amscd,typearea,hyperref}
\usepackage{amssymb}
\usepackage{a4wide}
\usepackage[mathscr]{eucal}
\usepackage{mathrsfs}
\usepackage{typearea}
\usepackage{charter}
\usepackage{pdfsync}

\usepackage{amscd,amsxtra,calc}
\usepackage{cmmib57}
\usepackage{url}

\usepackage[a4paper,width=16.2cm,top=3cm,bottom=3cm]{geometry}

\numberwithin{equation}{section}

\newtheorem{theorem}{Theorem}[section]

\newtheorem{proposition}[theorem]{Proposition}

\newtheorem{lemma}[theorem]{Lemma}
\newtheorem{remark}[theorem]{Remark}

\newcommand{\cali}[1]{\mathscr{#1}}

\newcommand{\dist}{\mathop{\mathrm{dist}}\nolimits}

\newcommand{\ddc}{dd^c}
\newcommand{\dc}{d^c}

\def\mC{\mathcal{C}}
\def\mD{\mathcal{D}}

\def\mV{\mathcal{V}}
\def\mW{\mathcal{W}}
\def\mU{\mathcal{U}}

\def\mR{\mathcal{R}}

\newcommand{\dbar}{\overline\partial}
\newcommand{\ddbar}{\partial\overline\partial}

\newcommand{\Cc}{\cali{C}}

\newcommand{\Fc}{\cali{F}}
\newcommand{\Gc}{\cali{G}}
\newcommand{\Hc}{\cali{H}}

\newcommand{\E}{\mathbb{E}}
\newcommand{\C}{\mathbb{C}}

\newcommand{\N}{\mathbb{N}}

\newcommand{\R}{\mathbb{R}}
\renewcommand\P{\mathbb{P}}

\title[Super-potentials, densities of currents and periodic points]{Super-potentials, densities of currents and number of\break periodic points for holomorphic maps}

\author{Tien-Cuong Dinh}
\address{Department of Mathematics, National University 
of Singapore, 10 Lower Kent Ridge Road, Singapore 119076.}
\email{matdtc@nus.edu.sg ; http://www.math.nus.edu.sg/$\sim$matdtc}
\author{Vi{\^e}t-Anh Nguy{\^e}n}
\address{Universit\'e de Lille 1, 
Laboratoire de math\'ematiques Paul Painlev\'e, 
CNRS U.M.R. 8524,  
59655 Villeneuve d'Ascq Cedex, 
France.}
\email{Viet-Anh.Nguyen@math.univ-lille1.fr} 
\author{Duc-Viet  Vu}
\address{Korea institute for advanced study,
 85 Hoegiro, Dongdaemun-gu, Seoul 02455, Republic of Korea}
\email{vuviet@kias.re.kr}

\dedicatory{Dedicated to  Professor  L\^e  Tu\^an   Hoa  on the occasion of his  sixtieth birthday}

\date{September 28, 2017}
\begin{document}
\begin{abstract}
We prove that if a positive closed current is bounded by another one  with bounded, continuous or H\"older continuous super-potentials, then it inherits the same property. 
There are two different methods to define wedge-products of positive closed currents of arbitrary bi-degree on compact K\"ahler manifolds using super-potentials and densities. When the first method applies, we show that the second method also applies and gives the same result. As an application, we obtain a sharp upper bound for the number of isolated periodic points of holomorphic maps on compact K\"ahler manifolds whose actions on cohomology are simple. A similar result still holds for a large class of holomorphic correspondences.
\end{abstract}
\maketitle

\medskip

\noindent
{\bf Classification AMS 2010}: 32U40, 32H50, 37F05.

\medskip


%\tableofcontents

%%%%%%%%%%%%%%%%%%%%%%%%%%%%%

%%%%%%%%%%%%%%%%%%%%%%%%%%%%%%%%%%
\section{Introduction} \label{S:introduction}

\maketitle

Positive closed currents are  a fundamental tool in complex analysis, algebraic
geometry, differential geometry and complex dynamics. Bi-degree
$(1, 1)$-currents and their intersections were intensively studied and have had many
applications while arbitrary bi-degree currents are much less understood, see Bedford-Taylor \cite{BedfordTaylor}, Demailly \cite{Demailly12, Demailly_survey} and Fornaess-Sibony \cite{FornaessSibony}.
In \cite[p.16]{Demailly92}, Demailly posed the problem to develop a theory of intersection
for positive closed currents of higher bi-degree.
 Sibony and the first author have  recently introduced,
in a series of articles \cite{DinhSibony09,DinhSibony10,DinhSibony12}, two   different  approaches 
 to    this problem  in the  context of compact K\"ahler manifolds. Several applications in complex dynamics and foliation theory have been obtained as well.
  
Let $X$ be a compact K\"ahler manifold of dimension $k$. 
Let $T$ and $S$ be  positive closed currents of bi-degree respective $(p,p)$  and $(q,q)$  on $X$.  
We will briefly recall the notions of super-potentials and densities of positive closed currents. We then describe the two different approaches to define  the  wedge-product (intersection) of $T$ and $S$.   Our first aim is to prove that 
when both methods apply we obtain the same current. This will allow us to unify the advantages of both methods. We also consider a domination principle for super-potentials and apply our study to bound the number of isolated periodic points of  holomorphic maps or correspondences.

The  notion of  super-potentials of positive closed currents were introduced in \cite{DinhSibony09,DinhSibony10}. We refer to these references for details. The starting point is that the
pluripotential theory is well developed for positive closed currents of bi-degree $(1, 1)$
thanks to the notion of plurisubharmonic (p.s.h. for short) functions. More precisely,
if $T$ is a positive closed $(1, 1)$-current, then we can write locally $T = \ddc u$ in the sense
of currents, where $u$ is a p.s.h. function, $\dc := {i\over 2\pi} (\partial - \dbar)$ and $\ddc={i\over\pi}\ddbar$. Since the function $u$ is everywhere defined,
if $u$ is integrable with respect to the trace measure of $S$ (we will say ``with respect to $S$" for short), then one can define $uS$ in the sense of currents and define $T\wedge S:=\ddc(uS)$. It is not difficult to check that the definition is independent of the choice of $u$ because $u$ is unique up to a pluriharmonic function.
When $T$ is of higher
bi-degree, we still can write locally $T=\ddc U$ but the potential $U$ does not satisfy the above properties of p.s.h. functions and one cannot consider their wedge-product in the same way.

Super-potentials are functions which play the role of quasi-potentials for positive closed
currents of arbitrary bi-degree.  For simplicity, we will not introduce this notion in full generality but limit ourselves in the necessary setting. Let $\mD_q(X)$ (or $\mD_q$ for short) denote the real vector space spanned by positive closed $(q,q)$-currents on $X$. 
Define the {\it $\ast$-norm} on this space by $\|R\|_\ast:=\min ( \|R^+\|+\|R^-\|)$, where $R^\pm$ are positive closed $(q,q)$-currents satisfying $R=R^+-R^-$ and $\|\ \|$ denotes the mass of a current. We consider this space of currents with the following {\it topology} : a sequence $(R_n)_{n\geq 0}$ in $\mD_q(X)$ converges in this space to $R$ if $R_n\to R$ weakly and if $\|R_n\|_\ast$ is bounded independently of $n$. On any $\ast$-bounded set of $\mD_q(X)$, this topology coincides with the classical weak topology for currents. 
It was shown in \cite{DinhSibony04} that the subspace $\widetilde\mD_q(X)$ of real closed smooth $(q,q)$-forms is dense in $\mD_q(X)$ for the considered topology, see Theorem \ref{t:reg} below.

Let $\mD_q^0(X)$ and $\widetilde\mD^0_q(X)$ denote the linear subspaces in $\mD_q(X)$ and $\widetilde\mD_q(X)$ respectively of currents whose cohomology classes in $H^{q,q}(X,\R)$ vanish. Their co-dimensions are equal to the dimension of $H^{q,q}(X,\R)$ which is finite. Fix a real smooth and closed $(p,p)$-form $\alpha$ in the same cohomology class with $T$ in $H^{p,p}(X,\R)$.
We will consider in this paper the super-potential of $T$ which is the real function $\mU_T$ on $\widetilde\mD^0_{k-p+1}$ defined by 
$$\mU_T(R):=\langle T-\alpha, U_R\rangle \qquad \text{for} \qquad R\in \widetilde\mD^0_{k-p+1},$$
where $U_R$ is any smooth form of bi-degree $(k-p,k-p)$ such that $\ddc U_R=R$. This form always exists because the cohomology class of $R$ vanishes. Note that since the cohomology class of $T-\alpha$ vanishes, we can write $T-\alpha=\ddc U_T$ for some current $U_T$. By Stokes theorem, we have  
$$\mU_T(R)=\langle \ddc U_T,U_R\rangle =\langle U_T,\ddc U_R\rangle =\langle U_T, R\rangle.$$
We deduce from this identities that $\mU_T(R)$ doesn't depend on the choice of $U_R$ and $U_T$. However,  $\mU_T$  depends on the  reference form $\alpha.$ Note also that
if $T$ is smooth, it is not necessary to take $R$ and $U_R$ smooth.

We will not consider other super-potentials of $T$. They are some affine extensions of $\mU_T$ to $\mD_{k-p+1}$ or its extensions to some subspaces. The following notions do not depend on the choice of super-potential nor on the reference form $\alpha$. We say that $T$ has a {\it bounded super-potential} if $\mU_T$ is bounded on each $\ast$-bounded subset of $ \widetilde\mD^0_{k-p+1}$. 
We say that $T$ has a {\it continuous super-potential} if $\mU_T$ can be extended to a continuous function on $\mD^0_{k-p+1}$. 
Recall that $\mD_p(X)$ is a metric space. If $\gamma>0$ is a constant, define for $R\in\mD_q(X)$
$$\|R\|_{-\gamma}:=\sup\big\{ |\langle R,\phi\rangle|, \ \phi \text{ is a test form of bi-degree } (k-q,k-q) \text{ with } \|\phi\|_{\mC^\gamma}\leq 1\big\}.$$ 
The topology induced by this distance coincides with the above-considered topology on each $\ast$-bounded set. By interpolation theory, in each $\ast$-bounded set, for all $\gamma\geq \gamma'>0$, there is a constant $c>0$ such that 
$\|\ \|_{-\gamma}\leq \|\ \|_{-\gamma'}\leq c (\|\ \|_{-\gamma})^{\gamma'/\gamma}$. We say that $T$ has a {\it H\"older continuous super-potential} if $\mU_T$ is continuous and H\"older continuous with respect to one of (or equivalently to all) norms $\|\ \|_{-\gamma}$ described above, on each $\ast$-bounded subset of  $\mD^0_{k-p+1}$.
Here is the domination principle for super-potentials that we mentioned above.

\begin{theorem}\label{t:main_1}
Let $T$ and $T'$ be positive closed $(p,p)$-currents on a compact K\"ahler manifold $X$. Assume that $T'\leq T$. 
If $T$ has a bounded, continuous or H\"older continuous  super-potential, then $T'$ satisfies the same property.
\end{theorem}

Consider now two positive closed currents $T$ and $S$ on $X$ of bi-degree $(p,p)$ and $(q,q)$ respectively. Assume that $p+q\leq k$ and that $T$ has a continuous super-potential. So $\mU_T$ is defined on whole $\mD^0_{k-p+1}$. We can define the wedge-product $T\wedge S$ by 
$$\langle T\wedge S,\phi\rangle := \langle \alpha\wedge S,\phi\rangle +\mU_T(S\wedge\ddc\phi)$$
for every smooth real test form $\phi$ of bi-degree $(k-p-q,k-p-q)$. Note that $S\wedge\ddc\phi$ belongs to $\mD^0_{k-p+1}$ because it is equal to $\ddc(S\wedge \phi)$. It is not difficult to check that $T\wedge S$ is equal to the usual wedge-product of $T$ and $S$ when one of them is smooth. The current $T\wedge S$ is positive and closed, see \cite{DinhSibony09, DinhSibony10,Vu}.

We discuss now the second definition of the wedge-product using the recent theory of densities for currents, see \cite{DinhSibony12} and Section \ref{s:densities} for details. The rough idea is that if we identify the diagonal $\Delta$ of $X\times X$ with $X$ in the canonical way, then the wedge-product between $T$ and $S$ will be the intersection between the current $T\otimes S$ with the current  of integration $[\Delta]$ on $\Delta$.  Let $\E$ denote the normal bundle to $\Delta$ in $X\times X$ and let $\pi:\E\to\Delta$ be the canonical projection. 
We  dilate roughly the current $T\otimes S$ near $\Delta$ using dilations in the normal directions to $\Delta$. This operation doesn't change the considered intersection with $[\Delta]$. 
When the dilation factor tends to infinity, the image of $T\otimes S$ admits limit values which are positive closed current of bi-degree $(p+q,p+q)$ in $\E$. These currents are called {\it the tangent currents} of $T \otimes S$ along $\Delta$. 

Assume there is only one tangent current and denote it by  $\mR$. We say that the {\it h-dimension of $\mR$ is less than $s$} if for any smooth $(s,s)$-form $\phi$ on $\Delta$ we have $\mR\wedge \pi^*(\phi)=0$. It is known that if the h-dimension of $\mR$ is less than $k-p-q+1$ then there is a unique positive closed current $R$ of bi-degree $(p+q,p+q)$ on $\Delta$ such that $\mR=\pi^*(R)$. We then define $T\curlywedge S:=R$ and say that   $T\curlywedge S$ is {\it well-defined}. The rough idea here is that $R$ is the intersection of $\mR$ with $[\Delta]$. 
The reader can consider the case where $T$ and $S$ are currents of integration on two submanifolds of $X$. Then the above condition on the h-dimension is equivalent to saying that the intersection of the two submanifolds is either empty or of minimal dimension $k-p-q$. 
This condition of the non-excess of dimension (which is a generic situation) is necessary to define the wedge-product as a current. 
Here is our second main theorem. It 
gives a new proof to the result by the third author saying that $T\wedge S$ is positive \cite{Vu}.

\begin{theorem}\label{t:main_2}
Let $X$ be a compact K\"ahler manifold of dimension $k$.
Let $T$ and $S$ be positive closed currents on $X$ of bi-degree $(p,p)$ and $(q,q)$ respectively. Assume that $p+q\leq k$ and that $T$ has a continuous super-potential. Then the wedge-product $T\curlywedge S$ is well-defined and we have 
 $T\wedge S= T\curlywedge S$. 
\end{theorem}

Our last result is an application of the above two theorems in complex   dynamics which is originally the motivation of this work.
A basic problem in dynamics is to study the distribution of periodic points when the period tends to infinity. 
We refer to  \cite{DinhSibony16b} for a  recent  survey on this  topic.
We will give here a sharp upper bound for the number of isolated periodic points of a holomorphic map on a compact K\"ahler manifold. Note that when the set of periodic points is of zero dimension then this result is an easy consequence of the classical Lefschetz's fixed point theorem. The existence of positive dimension sets of periodic points is the main difficulty in our problem, see Oguiso \cite{Oguiso} for some recent examples.

Let $f:X\to X$ be a holomorphic map which is {\it dominant}, i.e. surjective in the present setting. 
Recall that {\it the dynamical degree $d_s$ of order $s$} of $f$ is the spectral radius of the pull-back
operator $f^*$ acting on the Hodge cohomology group $H^{s,s} (X, \C)$ for $0 \leq s \leq  k$. It is known that $d_s$ itself is an eigenvalue of $f^*$ on $H^{s,s}(X,\C)$.
An inequality due to Khovanskii, Teissier and Gromov \cite{DinhNguyen06, Gromov} implies that the function
$s\mapsto  \log d_s$ is concave on $0 \leq  s\leq  k$. In particular, there are integers $p$ and $p'$ with $0 \leq p  \leq  p' \leq k$
such that
$$d_0 < \cdots < d_p = \cdots = d_{p'}  > \cdots > d_k .$$
We always have $d_0 =1$. The last dynamical degree
$d_k$  is also called {\it  the topological degree} of $f$ because it is equal to the  cardinality of $f^{-1}(x)$ for a  generic point $x$ in $X$. We call $d_p$ {\it the main dynamical  degree} of $f$.

In what follows, we assume that the action of $f^*$ on cohomology is {\it simple}, that is, $p=p'$ and 
$d_p$ is a simple eigenvalue and is the only eigenvalue of maximal  modulus of $f^*$ on $H^{p,p}(X,\C)$. In this case, it is known that $d_p$ is also a simple eigenvalue and is the only eigenvalue of maximal  modulus of $f^*$ on the full Hodge cohomology group $\oplus H^{*,*}(X,\C)$, see \cite[Prop. 5.8]{Dinh05} which holds for all holomorphic maps. Our  last main result is  the  following theorem, see also  our works with Truong \cite{DinhNguyenTruong15,DinhNguyenTruong16},  Diller-Favre \cite{Diller_Favre}, Favre \cite{Favre}, Iwasaki-Uehara \cite{IwasakiUehara}, Saito \cite{Saito}, Xie \cite{Xie} for similar and related statements in dimension 2 or for maps with dominant topological degree, i.e. the case where $p=p'=k$. 

\begin{theorem}\label{t:main_3}
 Let f be a dominant holomorphic endomorphism of a compact K\"ahler manifold $X$. Assume
that the action of the pull-back operator $ f^*ˆ—$ on cohomology is simple. Let $d$ be the main dynamical  degree of $f$ and let $P_n$ denote the number of isolated periodic  points of period $n$ of $f$ counted with multiplicity. 
Then  we have
$$P_n\leq  d^n +o(d^n) \quad \text{as} \quad  n\to\infty.$$
 \end{theorem}

 The plan of the paper is as follows. In Section \ref{s:sp}, basic facts on positive closed currents, super-potentials and the proof of Theorem \ref{t:main_1} will be given. Let $\Pi:\widehat{X\times X}\to X\times X$ be the blow-up of $X\times X$ along the diagonal $\Delta$ and let $\widehat\Delta:=\Pi^{-1}(\Delta)$ be the exceptional hypersurface. The manifold $\widehat{X\times X}$ has been used by Sibony and the first author to construct a good potential for currents, e.g. a solution $U_R$ for the above equation $\ddc U_R=R$, see e.g. \cite{DinhSibony09}. We can study the super-potentials using some integrals on $\widehat{X\times X}$. Some mass estimates of currents near $\widehat\Delta$ allow us to characterize currents with bounded, continuous or H\"older continuous super-potentials. We then deduce the desired regularity of $\mU_{T'}$ from the regularity of $\mU_T$ as stated in Theorem \ref{t:main_1}. 
 
In Section \ref{s:densities}, we will give basic facts from the theory of densities for currents and the proof of Theorem \ref{t:main_2}.
A key point is that the theory of densities allows us to see $T\curlywedge S$ using local holomorphic coordinates and smooth test forms. We use again 
$X\times X$ and  $\widehat{X\times X}$ in order to compare $T\curlywedge S$ with $T\wedge S$ and get the desired result. 

Finally, the proof of Theorem \ref{t:main_3} is given in Section \ref{s:dynamics} in a more general setting of correspondences on $X$. Let $\Gamma_n$ denote the graph of $f^n$ in $X\times X$ and $[\Gamma_n]$ the current of integration on $\Gamma_n$. 
Periodic points of period $n$ can be identified to the intersection between $\Gamma_n$ and $\Delta$. 
Using the above two main theorems, we show that the current $d^{-n}[\Gamma_n]$ converges to a positive closed current whose intersection with $\Delta$ has no "dimension excess". 
In some sense, although the intersection of $\Gamma_n$ with $\Delta$ may have positive dimension (this is the case when the set of periodic points has positive dimension), 
this dimension excess is negligible when $n$ tends to infinity. We are then in a situation similar to the one in our previous works  with  Truong 
\cite{DinhNguyenTruong15,DinhNguyenTruong16} and Theorem \ref{t:main_3} follows easily.

\medskip\noindent
{\bf Main notation.} Throughout the paper, $X$ is a compact K\"ahler manifold of dimension $k$. We also fix  a K\"ahler form $\omega$ on $X$. For simplicity, 
in each K\"ahler manifold endowed with a K\"ahler form, we use the Riemannian metric induced by the K\"ahler form. Let $\Pi:\widehat{X\times X}\to X\times X$ be the blow-up of $X\times X$ 
along the diagonal $\Delta$ and let $\widehat\Delta:=\Pi^{-1}(\Delta)$ be the exceptional hypersurface. Denote by $\pi_j$ the projections from $X\times X$ onto its factors for $j=1,2$ and define $\Pi_j:=\pi_j\circ\Pi$. The pairing $\langle\cdot,\cdot\rangle$ denotes the value of a current at a test form.
If $T$ is a positive or negative $(p,p)$-current on $X$, its mass is given by $\|T\|:=\langle T,\omega^{k-p}\rangle$ or  $\|T\|:=-\langle T,\omega^{k-p}\rangle$ respectively. The spaces $\mD_p,\mD_p^0, \widetilde\mD_p,\widetilde\mD_p^0$, their topology, the norms $\|\ \|_{-\gamma}$, $\|\ \|_*$, the super-potential $\mU_T$ and the wedge-products of currents are defined above. 
We also use $\{\ \}$ to denote the cohomology class of a current or an analytic set.  
 
\medskip

\noindent
{\bf Acknowledgments.} The  first author  was supported by Start-Up  Grant R-146-000-204-133 and Tier 1 Grant  R-146-000-248-114  from the National University of Singapore (NUS). The paper was partially prepared 
during the visits  of the authors at  the Freie Universit\"at Berlin, the
Vietnam  Institute for Advanced Study in Mathematics and 
at  the NUS.
They would like to express the gratitude to these institutions, the Alexander von Humboldt foundation and H\'el\`ene Esnault for their hospitality and  support.

%%%%%%%%%%%%%%%%%%%%%%%%%%%%%%%%%%%%%%%%%%%%%%%%%%%%%%%%%%%%%%%%%

\section{Domination principle for  super-potentials} \label{s:sp}

In this section, we will give some properties of positive closed currents, their super-potentials and then the proof of Theorem \ref{t:main_1}. 
The following regularization theorem for currents was obtained in \cite{DinhSibony04}. The last assertion in the statement is a consequence of the construction of $R_n^\pm$ in the above reference.

\begin{theorem} \label{t:reg}
Let $(X,\omega)$ be a compact K\"ahler manifold of dimension $k$. Let $R$ be a current in $\mD_q(X)$. Then there are positive closed $(q,q)$-currents $R^\pm$ and smooth positive closed $(q,q)$-forms $R_n^\pm$ in $\mD_q(X)$ such that 
$R=R^+-R^-$, $R_n^\pm$
are in the same cohomology class with $R^\pm$ in $H^{q,q}(X,\R)$ and $R^\pm_n$ converge weakly to $R^\pm$ as $n$ tends to infinity. Moreover, there is a constant $c>0$ independent of $R$ such that $\|R^\pm\|\leq c\|R\|_\ast$ and $\|R_n^\pm\|\leq c\|R\|_\ast$ for every $n$. In particular, $R_n:=R_n^+-R_n^-$ converges to $R$ in the considered 
topology on $\mD_q(X)$. If $R$ is a continuous form, then $R_n$ converges uniformly to $R$.
\end{theorem} 

Recall that we only use the super-potentials introduced in the introduction. We have the following basic  lemma.

\begin{lemma} \label{l:sp_pullback}
Let $\tau:\widetilde X\to X$ be a submersion between compact K\"ahler manifolds. Let $T$ be a positive closed $(p,p)$-current on $X$. 
If $T$ has a bounded, continuous or H\"older continuous super-potential, then the positive closed $(p,p)$-current 
$\tau^*(T)$ on $\widetilde X$ satisfies the same property.
\end{lemma}
\proof 
Let $l$ denote the dimension of $\widetilde X$. Let $\alpha$ be as in the introduction. Define $\mU:=\mU_T$, $\widetilde T:=\tau^*(T)$ and  $\widetilde\alpha:=\tau^*(\alpha)$.
So $\widetilde T$ is in the same cohomology class with $\widetilde\alpha$. We denote by $\widetilde \mU$ the super-potential of $\widetilde T$ defined as in the introduction using the reference form $\widetilde\alpha$. 

Recall that $\tau$ induces the push-forward operator $\tau_*$ on currents : if $\widetilde R$ is a current on $\widetilde X$ and if $\phi$ is a smooth test form of the right degree on $X$ then $\langle \tau_*(\widetilde R),\phi\rangle = \langle \widetilde R,\tau^*(\phi)\rangle$. 
This operator is continuous and is compatible with the operators $\partial$ and $\dbar$.  
So it is not difficult to see that $\tau_*$ defines a continuous map from $\mD_{l-p+1}^0(\widetilde X)$ to $\mD_{k-p+1}^0(X)$ and sends $\ast$-bounded sets to $\ast$-bounded sets. Moreover, given any positive real number $\gamma,$ the number  $\|\tau^*(\phi)\|_{\Cc^\gamma}$ is bounded by a constant times $\|\phi\|_{\Cc^\gamma}$. It follows that $\|\tau_*(\widetilde R)\|_{-\gamma}$ is bounded by a constant times $\|\widetilde R\|_{-\gamma}$, or in other words, 
the map $\tau_*: \mD_{l-p+1}^0(\widetilde X) \to \mD_{k-p+1}^0(X)$ is Lipschitz with respect to the norms $\|\ \|_{-\gamma}$. 

Consider any smooth form $\widetilde R$ in $\widetilde\mD_{l-p+1}^0(\widetilde X)$ and choose a smooth form $\widetilde U$ such that $\ddc \widetilde U=\widetilde R$. 
Define $R:=\tau_*(\widetilde R)$ and $U:=\tau_*(\widetilde U)$. We have $\ddc U=R$ and $R, U$ are smooth since $\tau$ is a submersion. In particular, $R$ belongs to $\widetilde\mD^0_{k-p+1}(X)$, or in other words, $\tau_*$ maps $\widetilde \mD_{l-p+1}^0(\widetilde X)$ to $\widetilde\mD_{k-p+1}^0(X)$. We also have 
$$\widetilde \mU(\widetilde R)=\langle \widetilde T-\widetilde\alpha, \widetilde U\rangle=\langle \tau^*(T-\alpha), \widetilde U\rangle
=\langle T-\alpha, U\rangle =\mU(R)$$
which implies the identity $\widetilde \mU=\mU\circ\tau_*$ on $\widetilde\mD_{l-p+1}^0(\widetilde X)$. 

We have seen that $\tau_*$ is bounded, continuous and Lipschitz. Therefore, if $\mU$ is bounded, continuous or H\"older continuous, $\widetilde\mU$ satisfies the same property. This completes the proof of the lemma.
\endproof

Recall that the projections $\Pi_j:\widehat{X\times X}\to X$ are submersions, see e.g. \cite{DinhSibony04}. We will apply the above lemma to these maps.
We have seen in the introduction that the definition of super-potential involves the solutions of the equation $\ddc U_R=R$. We will recall here the construction of kernel solving this solution and refer to \cite{DinhSibony10} for details. 

By Blanchard's theorem \cite{Blanchard}, $\widehat{X\times X}$ is a K\"ahler manifold. So we fix
a K\"ahler form $\widehat\Omega$ on $\widehat{X\times X}$ and we
can apply Hodge theory to this manifold. 
By K\"unneth's formula, the cohomology class $\{\Delta\}$ of $[\Delta]$ in $H^{k,k}(X \times X,\R)$ can be represented by a real smooth closed $(k,k)$-form $\alpha_{\Delta}$  which is a finite sum of forms of type $\pi_1^*(\phi_1) \wedge \pi_2^*(\phi_2)$ where $\phi_1$ and $\phi_2$ are closed smooth forms of suitable bi-degrees on $X$. By  \cite[1.2.1]{BGS} and \cite[Ex. 2.3.1]{DinhSibony10}, there is a real smooth closed $(k-1,k-1)$-form $\hat \eta$ on $\widehat{X\times X}$ such that $\hat\eta \wedge [\widehat\Delta]$ is cohomologous to $\Pi^*(\alpha_\Delta)$ and 
\begin{align} \label{e:Delta}
\Pi_{*}(\hat\eta \wedge [\widehat\Delta])=[\Delta].
\end{align}

Choose a real smooth closed $(1,1)$-form $\hat\beta$ on $\widehat{X\times X}$ which is cohomologous to $[\widehat\Delta]$. So we can write $[\widehat\Delta]-\hat\beta=\ddc \hat u$, where $\hat u$ is a quasi-p.s.h. function on $\widehat{X\times X}$. This equation implies that $ \hat u$ is smooth outside $\widehat\Delta$ and $\hat u-\log\dist(\cdot,\widehat\Delta)$ is a smooth function near $\widehat\Delta$. Subtracting from $\hat u$ a constant allows us to assume that $\hat u$ is negative.
Observe that since $\hat\beta\wedge\hat\eta$ is cohomologous to $\Pi^*(\alpha_\Delta)$, there is a real smooth $(k-1,k-1)$-form $\hat\beta'$ on $\widehat{X\times X}$ such that 
\begin{align} \label{e:beta}
\ddc \hat\beta'=\hat\beta\wedge\hat\eta -\Pi^*(\alpha_\Delta).
\end{align}
Adding to $\hat\beta'$ a constant times $\widehat\Omega^{k-1}$ allows us to assume that $\hat\beta'$ is positive. 
For $\epsilon>0$, denote by $\widehat\Delta_\epsilon$ the set of points in $\widehat{X\times X}$ with distance less than $\epsilon$ to $\widehat\Delta$. 
\begin{proposition} \label{p:kernel}
Let $R$ be a current in $\mD_q(X)$ such that $\|R\|_*\leq 1$. Then there is a constant $c>0$  independent of $R$ such that  the following properties hold for $0<\epsilon<1/2$.
\begin{itemize}
\item[(1)] $\hat u\Pi_2^*(R)$ is a current depending continuously on $R$ whose mass is bounded by $c$;
\item[(2)] The mass of $\hat u\Pi_2^*(R)$ on $\widehat\Delta_\epsilon$ is smaller than $c\epsilon^2|\log\epsilon|$;
\item[(3)] If $R$ is in  $\mD^0_q(X)$ and 
$U_R:=(\Pi_1)_*\big((\hat u\hat\eta+\hat\beta')\wedge\Pi_2^*(R)\big)$, then
$\ddc U_R=R$.
\end{itemize}
\end{proposition}

\noindent
{\bf Sketch of the proof.} This statement  was essentially 
obtained in \cite{DinhSibony10}. For the reader's convenience, we sketch here the proof. 
For the assertions (1) and (2), we can assume that $R$ is positive and we don't need to assume that it is closed. So it is enough to consider the case where $R$ is supported by a point $a\in X$ because all positive currents can be obtained as an average of such currents. The statement is now clear because $\Pi_2$ is a submersion and its fibers are transverse to the hypersurface $\widehat\Delta$. We used here the property of $\hat u$ described above and the fact that $\Pi_2^*(a)$ can be identified to the blow-up of $X$ at $a$ which varies continuously when $a$ varies.

Consider now the assertion (3). Combining the continuity obtained in the first assertion with Theorem \ref{t:reg}, it is enough to consider the case where $R$ is smooth. A direct computation using \eqref{e:beta} and the definition of $\hat u$ gives
$$\ddc U_R=(\Pi_1)_*\big([\widehat\Delta]\wedge\hat\eta\wedge\Pi_2^*(R)\big)- (\Pi_1)_*\big(\Pi^*(\alpha_\Delta)\wedge\Pi_2^*(R)\big).$$
Since the restriction of $\Pi_2^*(R)$ to $\widehat\Delta$ is equal to the one of $\Pi_1^*(R)$, using \eqref{e:Delta} and the identity $\Pi_1=\pi_1\circ\Pi$, we see that the first term in the RHS of the last equation is equal to
$$(\Pi_1)_*([\widehat\Delta\wedge\hat\eta])\wedge R= (\pi_1)_*[\Delta]\wedge R=R.$$

It remains to check that the second term vanishes. Using $\Pi_j=\pi_j\circ\Pi$, we see that this term satisfies
$$(\Pi_1)_*\big(\Pi^*(\alpha_\Delta)\wedge\Pi_2^*(R)\big)=(\pi_1)_*\Pi_*\big(\Pi^*(\alpha_\Delta\wedge\pi_2^*(R))\big)=(\pi_1)_*(\alpha_\Delta\wedge\pi_2^*(R)).$$
Now, if $\Phi$ is a smooth test form of the right bi-degree and if $(x,y)$ denotes the coordinates of points in $X\times X$, we have
$$\big\langle (\pi_1)_*(\alpha_\Delta\wedge\pi_2^*(R)),\Phi\big\rangle = \int_{X\times X} \Phi(x)\wedge \alpha_\Delta(x,y)\wedge R(y).$$
Since the cohomology class of $R$ vanishes, it is an exact form. Recall from the choice of  $\alpha_\Delta$ that it has  a nice  property of variable separation. Therefore, using Stokes and Fubini's theorems, we see that the last integral vanishes when we first integrate in variable $y$. 
This completes the proof of the proposition.
\hfill $\square$

\medskip

The following lemma allows us to compute the values of super-potentials.

\begin{lemma} \label{l:sp-comp}
Let $T,\alpha$ and $\mU_T$ be as above. Then for every $R$ in $\widetilde\mD^0_{k-p+1}(X)$, we have 
$$\mU_T(R)=  \int_{\widehat{X\times X}}   (\hat u\hat\eta+ \hat\beta') \wedge \Pi_1^*(T-\alpha)\wedge \Pi_2^*(R).$$
\end{lemma}
\proof 
Note that by Proposition \ref{p:kernel} applied to $T-\alpha, \Pi_1$ instead of $R,\Pi_2$, the integral in the lemma is meaningful. 
Let $U_R$ be as in Proposition \ref{p:kernel}. Observe that it is smooth since $R$ is smooth. By definition of super-potential in the introduction, we deduce from the definition of $U_R$ that
$$\mU_T(R)=\langle T-\alpha, U_R\rangle =\big\langle \Pi_1^*(T-\alpha), (\hat u\hat\eta+\hat\beta')\wedge \Pi_2^*(R)\big\rangle.$$
The lemma follows.
\endproof

The following proposition gives us a characterization of currents with bounded super-potentials.

\begin{proposition}\label{p:sp-bounded}
Let $T$ be a positive closed $(p,p)$-current on $X$. Then $T$ has a bounded super-potential if and only if there is a constant $c>0$ such that for every smooth positive closed $(k-p+1,k-p+1)$-form $R$ on $X$ with $\|R\|\leq 1$ we have
$$-\int_{\widehat{X\times X}} \hat u\, \widehat\Omega^{k-1} \wedge \Pi_1^*(T)\wedge \Pi_2^*(R)\leq c.$$
\end{proposition}
\proof
We first assume the last inequality and prove that $T$ has a bounded super-potential.
Let $R'$ be an arbitrary current in $\widetilde\mD_{k-p+1}^0(X)$ with $\|R'\|_*\leq 1$. 
We need to show that $|\mU_T(R')|$ is bounded by a constant. 
By Lemma \ref{l:sp-comp}, we have 
$$\mU_T(R')=  \int_{\widehat{X\times X}}   (\hat u\hat\eta+ \hat\beta') \wedge \Pi_1^*(T-\alpha)\wedge \Pi_2^*(R').$$

Observe that the $\ast$-norm of $\Pi_2^*(R')$ is bounded by a constant. By Proposition \ref{p:kernel} applied to $R'$ instead of $R$, the mass of  $\hat u\Pi_2^*(R')$ is bounded by a constant. We deduce that $|\mU_T(R')|$ is bounded by a constant if and only if 
$$\Big|\int_{\widehat{X\times X}}   \hat u\hat\eta \wedge \Pi_1^*(T)\wedge \Pi_2^*(R')\Big|$$
is bounded by a constant.
By the last assertion in Theorem \ref{t:reg}, it is enough to prove a similar property with a positive smooth closed form instead of $R'$.
So we can replace $R'$ by the form $R$ as in the proposition. Now it is clear that the inequality in the proposition implies the result because $\hat u$ is negative, $T$ and $R$ are positive, and we can bound $\hat\eta$ from above and below by constants times $\widehat\Omega^{k-1}$. Thus, $T$ has a bounded super-potential.

Assume now that $T$ has a bounded super-potential. We show that the inequality in the proposition holds for some constant $c$. 
Since $\|R\|\leq 1$, the cohomology class of $R$ is bounded and we can choose a real smooth closed form $\alpha_R$ in the cohomology class of $R$ whose $\Cc^0$-norm is bounded by a constant. Arguing as above, we only need to check that 
$$-\int_{\widehat{X\times X}} \hat u\, \widehat\Omega^{k-1} \wedge \Pi_1^*(T-\alpha)\wedge \Pi_2^*(R-\alpha_R)$$
is bounded from above by a constant.

Consider a smooth convex increasing function $\chi:\R\to\R$ vanishing on $(-\infty,-1]$ and equal to identity on $[1,+\infty)$. Define 
$\hat u_n:=\chi(\hat u+n)-n$. This function is smooth, equal to $\hat u$ on the set $\{\hat u\geq -n+1\}$ and decreases to $\hat u$ when $n$ tends to infinity. Moreover, $\|\ddc \hat u_n\|_*$ is bounded by a constant independent of $n$, see \cite[p.962]{DinhSibony04}. Therefore, in the last integral, we can replace $\hat u$ by $\hat u_n$. We only need here that the constants involving in our estimates  do not depend on $n$ nor on $R$.

Since the mass of a positive closed current depends only on its cohomology class, we easily deduce that the $\ast$-norm of 
$$\widetilde R':= \ddc \big(\hat u_n \widehat\Omega^{k-1} \wedge \Pi_2^*(R-\alpha_R)\big)$$
is bounded by a constant independent of $n$ and $R$. Observe also that the integral we consider is the value of the super-potential of $\Pi_1^*(T)$ at $\widetilde R'$. This value is bounded by a constant because by Lemma \ref{l:sp_pullback}, the current $\Pi_1^*(T)$ has a bounded super-potential. This completes the proof of the proposition.
\endproof

\noindent
{\bf Proof of Theorem \ref{t:main_1} for bounded super-potentials.} Assume that $T$ has a bounded super-potential. So $T$ satisfies the estimate in Proposition \ref{p:sp-bounded}. Since $\hat u$ is negative, $\widehat\Omega, T,R$ are positive and $T'\leq T$, the same property holds for $T'$ instead of $T$ with the same constant $c$.  Applying again Proposition \ref{p:sp-bounded} to $T'$ instead of $T$, we get that  $T'$ has a bounded super-potential.
\hfill $\square$

\medskip

\begin{remark} \label{r:main_1} \rm
Let $T, \alpha$ and $\mU_T$ be as above and assume that $T$ has a bounded super-potential. Let $\alpha'$ be a real smooth closed $(p,p)$-form. Then there is a constant $c>0$ such that for every positive closed $(p,p)$-current $T'$,
 cohomologous to $\alpha'$, with $T'\leq T$, we have 
$|\mU_{T'}(R)|\leq c \|R\|_*$ for $R\in\widetilde\mD_{k-p+1}(X)$. Here,  
 $\mU_{T'}$ denotes the super-potential of $T'$ associated with the reference form $\alpha'$. 
 \end{remark}

We give now a criterion to check if a current has a continuous super-potential. Let $T$ be a positive closed $(p,p)$-current as above. Recall that
for $\epsilon>0$, $\widehat\Delta_\epsilon$ denotes the set of points in $\widehat{X\times X}$ with distance less than $\epsilon$ to $\widehat\Delta$. 
Consider the following quantity
$$\vartheta_T(\epsilon):=\sup_R \int_{\widehat\Delta_\epsilon} -\hat u\widehat\Omega^{k-1}\wedge \Pi_1^*(T)\wedge \Pi_2^*(R),$$
where  the supremum is taken over all smooth positive closed forms $R$ on $X$, of bi-degree $(k-p+1,k-p+1)$, such that $\|R\|\leq 1$.

\begin{proposition} \label{p:sp-continuous}
Let $T$ be a positive closed $(p,p)$-current on $X$. Then $T$ has a continuous super-potential if and only if $\vartheta_T(\epsilon)$ tends to $0$ as $\epsilon$ tends to $0$. 
\end{proposition}
\proof
Assume that $T$ has a continuous super-potential. We prove that $\vartheta_T(\epsilon)$ tends to $0$ as $\epsilon$ tends to $0$. 
Fix an integer $n_0$ large enough so that $|\hat u-\log\dist(\cdot,\widehat\Delta)|\leq n_0$. 
Define $\epsilon_n:=e^{-2n-3n_0}$ and $\vartheta_n:=\vartheta_T(\epsilon_n)$. 
Since the function $\epsilon\mapsto \vartheta_T(\epsilon)$ is increasing, it is enough to show that $\vartheta_n$ tends to 0 as $n$ goes to infinity.  Let $\hat u_n$ be as in  the proof of Proposition \ref{p:sp-bounded}. Observe that on the set $\widehat\Delta_{\epsilon_n}$ we have
$$\hat u_n-\hat u = - n - \hat u \geq -n-\log\dist(\cdot,\widehat\Delta)-n_0\geq {1\over 2} (-\log\dist(\cdot,\widehat\Delta)+n_0)\geq -{1\over 2}\hat u.$$
Therefore, by the definition of $\vartheta_T(\epsilon)$, it is enough to show that 
$$\int_{\widehat{X\times X}} (\hat u_n-\hat u) \widehat\Omega^{k-1}\wedge \Pi_1^*(T)\wedge \Pi_2^*(R)$$
converges to 0 uniformly in $R$. 

We will replace $\hat u$ in the last integral by $\hat u_{n+m}$. It is enough to check that the obtained integral tends to 0 uniformly in $R$ and $m,n$. 
By Proposition \ref{p:kernel}, applied to $T,\Pi_1$ instead of $R,\Pi_2$, this property is true if we replace $T$ by any smooth form. Therefore, we only have to check that 
$$\int_{\widehat{X\times X}} (\hat u_n-\hat u_{n+m}) \widehat\Omega^{k-1}\wedge \Pi_1^*(T-\alpha)\wedge \Pi_2^*(R)$$
converges to 0 uniformly in $R$ and $m,n$. 
Denote by $\widetilde\mU$ the super-potential of $\Pi_1^*(T)$ defined as in the introduction using the reference form $\Pi_1^*(\alpha)$. The last expression is just the value of $\widetilde\mU$ at the current 
$$\widetilde R_{n,m} :=\ddc (\hat u_n-\hat u_{n+m}) \wedge \widehat\Omega^{k-1}\wedge \Pi_2^*(R)=\ddc \big[(\hat u_n-\hat u_{n+m})\widehat\Omega^{k-1}\wedge \Pi_2^*(R)\big].$$
As in the proof of Proposition \ref{p:sp-bounded}, this current has a $\ast$-norm bounded by a constant. Therefore, by Proposition \ref{p:kernel}, it converges (weakly) to 0 in $\mD_{2k-p+1}(\widehat{X\times X})$ as $n$ goes to infinity, uniformly in $R$ and $m$. Thus, $\widetilde\mU(\widetilde R_{n,m})$ tends to $\widetilde\mU(0)=0$ because by Lemma \ref{l:sp_pullback}, $\widetilde\mU$ is continuous. This completes the first part of the proof.

Assume now that $\vartheta_T(\epsilon)$ tends to $0$ as $\epsilon$ tends to $0$. We need to show that $T$ has a continuous super-potential. We first prove that it has a bounded super-potential. 
Consider a  positive closed form $R$ of bi-degree $(k-p+1,k-p+1)$  such that $\|R\|\leq 1$. 
By Proposition \ref{p:sp-bounded}, it is enough to show that 
$$\int_{\widehat{X\times X}} -\hat u\widehat\Omega^{k-1} \wedge \Pi_1^*(T)\wedge \Pi_2^*(R)$$
is bounded by a constant independent of $R$. 

Recall again that the mass of a positive closed current depends only on its cohomology class. Therefore, the mass of $ \Pi_1^*(T)\wedge \Pi_2^*(R)$ is bounded by a constant independent of $R$. 
Fix an $\epsilon_0>0$ such that $\vartheta_T(\epsilon_0)$ is finite.  The considered integral is equal to
$$-\Big[\int_{\widehat{X\times X}\setminus \widehat\Delta_{\epsilon_0}} +\int_{\widehat\Delta_{\epsilon_0}}  \Big]
 \hat u\widehat\Omega^{k-1} \wedge \Pi_1^*(T)\wedge \Pi_2^*(R).$$
The first integral is bounded because $\hat u$ is bounded on $\widehat{X\times X}\setminus\widehat\Delta_{\epsilon_0}$. The second one is bounded by $\vartheta_T(\epsilon_0)$. Thus, $T$ has a bounded super-potential.

We deduce from the boundedness of super-potential that if $R_l\in\widetilde\mD_{k-p+1}^0(X)$ is a $\ast$-bounded sequence, we can extract a subsequence $R_{l_j}$ such that $\mU_T(R_{l_j})$ admits a limit when $j$ tends to infinity. In order to show that $\mU_T$ is continuous, it is enough to show that the limit is independent of the choice of the subsequence. By linearity, we only have to show that if $R_l\in\widetilde \mD_{k-p+1}^0(X)$ converges weakly to 0 and $\|R_l\|_*\leq 1$ then $\mU_T(R_l)$ tends to 0 as $l$ tends to infinity. 

Fix an arbitrary positive number $\delta$ and then choose a number $0<\epsilon_0<\delta$ small enough such that  $\vartheta_T(4\epsilon_0)\leq \delta$. Let $0\leq \rho\leq 1$ be a smooth function with compact support in $\widehat\Delta_{2\epsilon_0}$ which is equal to 1 on $\widehat\Delta_{\epsilon_0}$.
By Lemma \ref{l:sp-comp}, we have 
$$\mU_T(R_l)=  \int_{\widehat{X\times X}}  \rho (\hat u\hat\eta+ \hat\beta') \wedge \Pi_1^*(T-\alpha)\wedge \Pi_2^*(R_l)
+ \int_{\widehat{X\times X}}   (1-\rho) (\hat u\hat\eta+ \hat\beta') \wedge \Pi_1^*(T-\alpha)\wedge \Pi_2^*(R_l).$$
For the first term in the last sum, we use  the definition of $\vartheta_T(2\epsilon_0)$ in order to bound the part involving $T$ and
 the second assertion of  Proposition \ref{p:kernel} to bound the part involving $\alpha$. We see that this term is smaller than $\delta$ times a constant independent of $l$. Since $\delta$ is arbitrary, it suffices to check that the second term tends to 0 as $l$ tends to infinity.

Since $(1-\rho) (\hat u\hat\eta+ \hat\beta')$ is smooth and vanishes near $\widehat\Delta$, it is equal to $\Pi^*(\Phi)$ for some smooth form $\Phi$ on $X\times X$ which vanishes near $\Delta$. It follows that the considered term is equal to 
$\big\langle (T-\alpha)\otimes R_l,\Phi\big\rangle$. Recall that the tensor product of two currents depends continuously on these currents. Since $R_l$ converges weakly to 0, the last pairing converges to 0. This completes the proof of the proposition.
\endproof

\noindent
{\bf Proof of Theorem \ref{t:main_1} for continuous super-potentials.} 
Assume that $T$ has a continuous super-potential. By Proposition \ref{p:sp-continuous}, $\vartheta_T(\epsilon)$ tends to $0$ as $\epsilon$ tends to $0$. 
It is clear that $\vartheta_{T'}(\epsilon)\leq \vartheta_T(\epsilon)$. So we also have that $\vartheta_{T'}(\epsilon)$ tends to $0$.
By applying again Proposition \ref{p:sp-continuous}, we obtain that $T'$ has a continuous super-potential.
\hfill $\square$

\medskip

The following result gives us a characterization of currents with H\"older continuous super-potentials.

\begin{proposition} \label{p:sp-Holder}
Let $T$ be a positive closed $(p,p)$-current on $X$. Then $T$ has a H\"older continuous super-potential if and only if $\vartheta_T(\epsilon)=O(\epsilon^\nu)$ as $\epsilon\to 0$ for some constant $\nu>0$.
\end{proposition}
\proof
Assume that $T$ has a H\"older continuous super-potential. We show that $\vartheta_T(\epsilon)=O(\epsilon^\nu)$ as $\epsilon\to 0$ for some constant $\nu>0$. As in the proof of Proposition \ref{p:sp-continuous}, it is enough to show that $|\mU(\widetilde R_{n,m})| \leq ce^{-\nu n}$ for some positive constants $c$ and $\nu$ independent of $R$ and $m$. By Lemma \ref{l:sp_pullback}, $\Pi_1^*(T)$ has a H\"older continuous super-potential. So it is enough to check that $\|\widetilde R_{n,m}\|_{-2}$ is bounded by a constant times $e^{-\nu n}$ for some positive constant $\nu$. For this purpose, by Theorem \ref{t:reg}, we can assume, without loss of generality, that $R$ is a positive closed form of mass 1. 

Let $\widetilde\Phi$ be a smooth test form of the right bi-degree in $\widehat{X\times X}$ with $\|\widetilde\Phi\|_{\Cc^2}\leq 1$. We have to show that $\big|\big\langle \widetilde R_{n,m}, \widetilde\Phi\big\rangle \big|$ is bounded by a constant times $e^{-n}$. We have
$$\big|\big\langle \widetilde R_{n,m}, \widetilde\Phi\big\rangle \big| = \big|\big\langle (\hat u_n-\hat u_{n+m}) \wedge \widehat\Omega^{k-1}\wedge \Pi_2^*(R), \ddc \widetilde\Phi\big\rangle \big|.$$
Observe that the form $\ddc\widetilde\Phi$ is bounded, $\widehat\Omega^{k-1}\wedge \Pi_2^*(R)$ is positive and $0\leq \hat u_n- \hat u_{n+m}\leq -\hat u$. Moreover, $\hat u_n- \hat u_{n+m}$ vanishes outside $\widehat\Delta_{e^{-n+n_0+1}}$ since both $\hat u_n$ and $\hat u_{n+m}$ are equal to $\hat u$ outside this set. Therefore,
the last pairing (which involves $\ddc\widetilde\Phi$) is bounded by a constant times the mass of $-\hat u  \Pi_2^*(R)$ in $\widehat\Delta_{e^{-n+n_0+1}}$. The second assertion in Proposition \ref{p:kernel} implies the desired estimate.

Assume now that $\vartheta_T(\epsilon)=O(\epsilon^\nu)$ as $\epsilon\to 0$ for some constant $\nu>0$. We prove that $T$ has a H\"older continuous super-potential. It follows from Theorem \ref{t:main_1} for continuous super-potentials that $T$ has a continuous super-potential. Let $R$ be a smooth form in $\widetilde\mD^0_{k-p+1}(X)$ with $\|R\|_*\leq 1$. We need to show that $|\mU_T(R)|\leq c'(\|R\|_{-2})^{\nu'}$ for some positive constants $c'$ and $\nu'$ independent of $R$. 
It is enough to consider the case where $\|R\|_{-2}$ is small. Define $\epsilon_0:=(\|R\|_{-2})^{1/(25k)}$. As in the end of 
Proposition \ref{p:sp-continuous}, consider a smooth function $0\leq \rho\leq 1$ with compact support in $\widehat\Delta_{2\epsilon_0}$ which is equal to 1 on $\widehat\Delta_{\epsilon_0}$.
We have 
$$\mU_T(R)=  \int_{\widehat{X\times X}}  \rho (\hat u\hat\eta+ \hat\beta') \wedge \Pi_1^*(T-\alpha)\wedge \Pi_2^*(R)
+ \int_{\widehat{X\times X}}   (1-\rho) (\hat u\hat\eta+ \hat\beta') \wedge \Pi_1^*(T-\alpha)\wedge \Pi_2^*(R).$$
The first term is bounded by a constant times $\vartheta(4\epsilon_0)+\epsilon_0$. Therefore, this term is bounded by $c'\|R\|_{-2}^{\nu'}$ for some positive constants $c'$ and $\nu'$. 

It remains to bound the second term. As in Proposition \ref{p:sp-continuous}, this term is equal to 
$$\big\langle (T-\alpha)\otimes R,\Phi\big\rangle=\big\langle T-\alpha, (\pi_1)_*(\pi_2^*(R) \wedge\Phi)\big\rangle.$$  
Observe that we can choose $\rho$ such that $\|\rho\|_{\Cc^2}$ bounded by a constant times $\epsilon_0^{-2}$. 
 The $\Cc^2$ norm of $\hat u$ is bounded by a constant times $\epsilon_0^{-2}$ on the support of $\Pi^*(\Phi)$ because this support is outside $\widehat\Delta_{\epsilon_0}$.
It follows that the $\Cc^2$ norm of $\Pi^*(\Phi)$ is bounded by $\epsilon_0^{-4}$.
Observe also that the $\Cc^3$ norm of the map $\Pi^{-1}$ on the support of $\Phi$ is bounded by a constant times $\epsilon_0^{-4}$. This can be seen using local coordinates in the construction of the blow-up of manifolds. We then deduce that the $\Cc^2$ norm of $\Phi$, which is equal to $(\Pi^{-1})^*\Pi^*(\Phi)$, is bounded by $\epsilon_0^{-24k}$. To see the last point, we can use a local real coordinate system $(t_1,\ldots,t_{4k})$ in $\widehat {X\times X}$. It is enough to observe that $\Cc^2$ norm of    $(\Pi^{-1})^*(dt_j)$ is $O(\epsilon_0^{-4})$ for each $1\leq j\leq 4k$ and if $g(t)$ is a coefficient of $\Pi^*(\Phi)$, the $\Cc^2$ norm of $g$ is $O(\epsilon_0^{-4})$ and the one of $(\Pi^{-1})^*(g)$ is $O(\epsilon_0^{-8})$ which is bounded by $O(\epsilon_0^{-8k})$.

The form $(\pi_1)_*(\pi_2^*(R) \wedge\Phi)$ is obtained from $\pi_2^*(R) \wedge\Phi$ by integrating on the fibers of $\pi_1$. Therefore, its $\Cc^0$ norm is bounded by a constant times $\|\Phi\|_{\Cc^2}\|R\|_{-2}$. By the definition of $\epsilon_0$, the last expression is bounded by a constant times $\epsilon_0$. Finally, since the mass of $T-\alpha$ is finite, we deduce that $\big|\big\langle T-\alpha, (\pi_1)_*(\pi_2^*(R) \wedge\Phi)\big\rangle\big|$ is bounded by a constant times $\epsilon_0$. This completes the proof of the proposition.
\endproof

\noindent
{\bf Proof of Theorem \ref{t:main_1} for H\"older continuous super-potentials.} 
Assume that $T$ has a H\"older continuous super-potential. Since $\vartheta_{T'}(\epsilon)\leq \vartheta_T(\epsilon)$, it is clear from Proposition \ref{p:sp-Holder} that $T'$ satisfies the same property. This ends the proof of Theorem \ref{t:main_1} for all cases.
\hfill $\square$

%%%%%%%%%%%%%%%%%%%%%%%%%%%%%%%%%%%%%%%%%%%%%%%%%%%%%%%%%%%%%%%%%
\section{Super-potentials versus  densities}\label{s:densities}

In this section, we first recall some basic facts from the theory of densities for currents and refer the reader to \cite{DinhSibony12} for details. For simplicity, we will restrict ourselves in the setting of Theorem \ref{t:main_2}. The proof of this theorem will be given at the end of the section.

Let $\pi:\E\to\Delta$ denote the normal vector bundle to $\Delta$ in $X\times X$. We identify $\Delta$ with the zero section of $\E$. For $\lambda\in\C^*$, let $A_\lambda:\E\to \E$ be the multiplication by $\lambda$ on the fibers of $\pi$. The diagonal $\Delta$ is invariant under the action of $A_\lambda$.
Consider a diffeomorphism $\tau$ from a neighborhood of $\Delta$ in $X\times X$ to a neighborhood of $\Delta$ in $\E$ such that the restriction of $\tau$ to $\Delta$ is the identity map. 
We see that 
the normal  bundle to $\Delta$ in $\E$    is  canonically  isomorphic  to  $\E$ and the differential $d\tau$ of $\tau$ preserves the tangent bundle of $\Delta.$ Consequently, $d\tau$ induces a real endomorphism of $\E$. There exists a map $\tau$ such that this endomorphism of $\E$ is the identity map, see \cite[Lemma 4.2]{DinhSibony12}. We say that $\tau$ is {\it admissible}. From now on, we fix an admissible map $\tau$ as above. In general, it is not a holomorphic map. 

Let $T$ and $S$ be as in Theorem \ref{t:main_2}. Consider the following family of closed currents of degree $(2p +2q)$ indexed by $\lambda\in\C^*$
$$\mR_\lambda:=(A_\lambda)_*\tau_*(T\otimes S).$$
Since $\tau$ may not be holomorphic, the current $\mR_\lambda$  may not be of bi-degree $(p+q,p+q)$ and we cannot talk about its positivity. However, for any sequence $(\lambda_n)_{n\geq1}$ in $\C^*$ converging to infinity, there is a subsequence $(\lambda_{n_j})_{j\geq 1}$ such that $\mR_{\lambda_{n_j}}$ converges to some positive closed current $\mR$ of bi-degree $(p+q,p+q)$ in $\E$, as $j$ tends to infinity. We say that $\mR$ is {\it a tangent current} of $T\otimes S$ along $\Delta$. It may depend on the sequence $\lambda_{n_j}$ but it is independent of the choice of $\tau$. Recall that $\mR$ is invariant under the action of $A_\lambda$ for every $\lambda\in\C^*$.

Tangent currents can be seen using local holomorphic coordinates near $\Delta$. We will introduce here coordinates which are suitable for the proof of Theorem \ref{t:main_2}. Let $z=(z_1,\ldots,z_k)$ denote a local holomorphic coordinate system on a local chart  $U$ of $X$. We consider the natural coordinate system $(z,w)$ on $U\times U$ with $w=(w_1,\ldots,w_2)$, a copy of $z$, such that $\Delta_U:=\Delta\cap (U\times U)$ is given by the equation $z=w$. We will use the coordinates $(z',w')$ for a small neighborhood $W$ of $\Delta_U$ with $z':=z$ and $w':=w-z$. So $\Delta_U$ is given by the equation $w'=0$. 
The restriction of $\E$ to $\Delta_U$ can be identified with $\Delta_U\times \C^k$. In this setting, the projection $\pi$ is just the map $(z',w')\mapsto (z',0)$. The dilation $A_\lambda$ is the map $(z',w')\mapsto (z',\lambda w')$. 

For simplicity, we also identify $W$ with an open subset of $\Delta_U\times \C^k$. With all these notations, the current $\mR$ above satisfies
$$\mR=\lim_{j\to\infty} (A_{\lambda_{n_j}})_*(T\otimes S) \quad \text{on} \quad \Delta_U\times \C^k$$
or equivalently, for any real smooth form $\Phi$ of bi-degree $(2k-p-q,2k-p-q)$ with compact support in $\Delta_U\times \C^k$, we have 
\begin{equation} \label{e:tangent}
\langle \mR, \Phi\rangle =\lim_{j\to\infty} \big\langle (A_{\lambda_{n_j}})_*(T\otimes S), \Phi\big\rangle.
\end{equation}
In what follows, we only need to consider $\lambda$ such that $|\lambda|\geq 1$. 

Let $i:X\to \Delta$ be the natural map $x\mapsto (x,x)$ sending $X$ to $\Delta$. 
In order to prove Theorem \ref{t:main_2}, it is enough to show that $\mR=\pi^*(i_*(T\wedge S))$, where $T\wedge S$ is defined using the super-potential of $T$ as in the introduction. 
So we need to study the pairing in the RHS of the identity \eqref{e:tangent}. 
We first consider the case where $\Phi$ has compact support in $W$ which is an open subset of $X\times X$.
We have (see also Lemma \ref{l:Psi-lambda} below)
\begin{equation} \label{e:Psi-lambda}
 \big\langle (A_\lambda)_*(T\otimes S), \Phi\big\rangle= \big\langle T\otimes S, A_\lambda^*(\Phi)\big\rangle=\big\langle T,\Psi_\lambda\big\rangle,
\end{equation}
where 
$$\Psi_\lambda:= (\pi_1)_*\big[A_\lambda^*(\Phi)\wedge \pi_2^*(S)\big].$$

\begin{lemma} \label{l:Psi-lambda}
The current $\Psi_\lambda$ is smooth, depending continuously on $S$, and the identity \eqref{e:Psi-lambda} holds. Moreover,
for $|\lambda|\geq 1$, 
the $L^1$-norm of $\Psi_\lambda$ is bounded by a constant independent of $\lambda$ and also of $S$ if we assume that the mass of $S$ is bounded by a fixed constant.
\end{lemma}
\proof
We only consider $\lambda$ such that $|\lambda|\geq 1$. 
We use the coordinates $z$ for $\Psi_\lambda$ and $(z',w')$ with $z'=z, w'=w-z$ for $\Phi$. We have 
$$\Psi_\lambda(z)=\int_{w\in X} \Phi\big(z,\lambda(w-z)\big) \wedge S(w).$$
Now, it is clear that $\Psi_\lambda$ is smooth and depends continuously on $S$.
The identity \eqref{e:Psi-lambda} clearly holds when $T$ is smooth. Since the tensor product of currents depend continuously on these currents, the identity holds for all $T$ and $S$. We don't need here that $T$ and $S$ are positive and closed.

For the rest of the lemma,  we only need to assume that $S$ is a positive current of bounded mass, not necessarily closed. 
We can assume that $S$ has support in a point $a\in X$ since it can be obtained as an average of such currents.  So $\pi_2^*(S)$ is supported by $X\times \{a\}$. Observe that the intersection of the support of $A_\lambda^*(\Phi)$
with $X\times \{a\}$ has $2k$-dimensional volume bounded by a constant times $|\lambda|^{-2k}$. Moreover, since the degree of $\Phi$ in $w'$ is at most $2k$, the coefficients of  $\Phi\big(z,\lambda(w-z)\big)$ is bounded by a constant times $|\lambda|^{2k}$. It is now clear that the $L^1$-norm of $\Psi_\lambda$ is bounded by a constant. 
\endproof

\begin{lemma} \label{l:Psi-lambda-lim}
We have that  $\Psi_\lambda$ converges weakly to $i^*(\pi_*(\Phi))\wedge S$ as $\lambda$ tends to infinity.
\end{lemma}
\proof
From the above study of $L^1$-norm of $\Psi_\lambda$ in the proof of Lemma \ref{l:Psi-lambda}, we see that if the degree of $\Phi$ in $w'$ is not maximal, then the $L^1$-norm of $\Psi_\lambda$ tends to 0. It is also clear that $\pi_*(\Phi)$ vanishes and the lemma holds in this case. So we can remove from $\Phi$ all terms which are not of maximal degree in $w'$. Without loss of generality, assume that 
$$\Phi(z',w')=\phi(z') \wedge h(z',w') (dw_1'\wedge d\overline w_1') \wedge \ldots \wedge  (dw_k'\wedge d\overline w_k'),$$
where $\phi$ is a smooth form of bi-degree $(k-p-q,k-p-q)$  on $U$ and $h(z',w')$ is a smooth function with compact support in $W$.

We will prove the lemma for any positive current $S$, not necessarily closed. 
Such a current $S$ can be written as an average of positive currents of the irreducible form
$$\delta_a \otimes (-iv_1\otimes \overline v_1)\otimes \cdots \otimes (-i v_{k-q}\otimes \overline v_{k-q}),$$
where $a$ in a point in $X$ and $v_j$ is a holomorphic tangent vector of $X$ at $a$. We have seen that the $L^1$-norm of $\Psi_\lambda$ is bounded by a constant independent of $\lambda$ and $S$. Therefore, by Lebesgue's convergence theorem, we can assume that $S$ has the above irreducible form. 

Note that if $a$ is outside $U$ then both $\Psi_\lambda$ and $i^*(\pi_*(\Phi))\wedge S$ vanish when $\lambda$ is large enough, and the lemma is clear. Consider the case where $a$ is in $U$. Using a linear change of coordinates, we can assume for simplicity that $z=0$ at $a$ and $v_j=\partial/\partial z_j$, that is, for any test form $\psi$ of bi-degree $(k-q,k-q)$ on $X$, we have $\langle S,\psi\rangle =(-i)^{k-q}g(0)$ if $g(z)$ denotes the coefficient of $dz_1\wedge d\overline z_1\wedge \ldots \wedge dz_{k-q}\wedge d\overline z_{k-q}$ in $\psi$.
Let $\Theta$ be a smooth $(p,p)$-form on $X$. We need to show that 
\begin{equation} \label{e:Psi-lambda-lim-main}
\lim_{\lambda\to\infty} \langle \Psi_\lambda,\Theta\rangle = \big\langle i^*(\pi_*(\Phi))\wedge S,\Theta\big\rangle.
\end{equation}

We first compute the limit in the last identity. By definition of $\Psi_\lambda$, we have 
$$\langle \Psi_\lambda,\Theta\rangle = \big\langle S, (\pi_2)_*\big(A_\lambda^*(\Phi)\wedge \pi_1^*(\Theta)\big)\big\rangle 
=(-i)^{k-q}g_\lambda(0),$$
where $g_\lambda(w)$ is the coefficient of $dw_1\wedge \ldots \wedge  d\overline w_{k-q}$ in $(\pi_2)_*\big(A_\lambda^*(\Phi)\wedge \pi_1^*(\Theta)\big)$.
In order to compute this coefficient, observe that 
$$A_\lambda^*(\Phi)\wedge \pi_1^*(\Theta)= |\lambda|^{2k} \phi(z) \wedge h(z, \lambda(w-z))
(dw_1'\wedge d\overline w_1') \wedge \ldots \wedge  (dw_k'\wedge d\overline w_k')\wedge \Theta(z).$$
Note that the image of the last form by $(\pi_2)_*$ is obtained by integrating it in variable $z$. Recall that  we only need to compute the coefficient of  $dw_1\wedge \ldots \wedge d\overline w_{k-q}$ at the point $w=0$. Since $w'=w-z$, we obtain
$$g_\lambda(0)=|\lambda|^{2k}\int_z \phi(z) \wedge h(z, -\lambda z) \ dz_{k-q+1}\wedge\ldots\wedge d\overline z_k\wedge \Theta(z).$$
Using the change of coordinates $z\mapsto -\lambda^{-1}z$ gives
$$g_\lambda(0)=\int_z |\lambda|^{2(k-p-q)}\phi(-\lambda^{-1}z) \wedge h(-\lambda^{-1}z,  z) \ dz_{k-q+1}\wedge \ldots \wedge d\overline z_k\wedge |\lambda|^{2p}\Theta(-\lambda^{-1}z).$$
It follows that
\begin{eqnarray*}
\lim_{\lambda\to\infty}  \langle \Psi_\lambda,\Theta\rangle & = & \lim_{\lambda\to\infty} (-i)^{k-q} g_\lambda(0)= (-i)^{k-q}\int_z \phi(0) \wedge h(0, z) \ dz_{k-q+1}\wedge \ldots \wedge d\overline z_k\wedge \Theta(0) \\
& = & (-i)^{k-q} \theta \int_z h(0, z) \ dz_1\wedge \ldots \wedge d\overline z_k,
\end{eqnarray*}
where $\theta$ is the coefficient of $dz_1\wedge\ldots \wedge d\overline z_{k-q}$ in $\phi(0)\wedge \Theta(0)$.

Consider now the RHS of \eqref{e:Psi-lambda-lim-main}. It is equal to $\big\langle S, i^*(\pi_*(\Phi))\wedge\Theta\big\rangle$. 
So it is equal to $(-i)^{k-q}l(0)$, where $l(z)$ is the coefficient of $dz_1 \wedge \ldots \wedge  d\overline z_{k-q}$ in
$ i^*(\pi_*(\Phi))\wedge\Theta$. From the above expression of $\Phi$ in coordinates $(z',w')$, we obtain that 
$$i^*(\pi_*(\Phi))= \phi(z) \int_{w'} h(z, w')dw_1'\wedge \ldots \wedge  d\overline w_k'.$$
It follows that 
$$\big\langle S, i^*(\pi_*(\Phi))\wedge\Theta\big\rangle=(-i)^{k-q}l(0)=(-i)^{k-q}\theta  \int_{w'} h(0, w')dw_1'\wedge \ldots \wedge  d\overline w_k',$$
which implies the desired identity \eqref{e:Psi-lambda-lim-main}.
\endproof

\begin{lemma} \label{l:ddc-Psi}
The $\ast$-norm $\|\ddc\Psi_\lambda\|_*$  of $\ddc\Psi_\lambda$ is bounded by a constant independent of $\lambda$ for $|\lambda|\geq 1$. In particular, $\ddc \Psi_\lambda$ converges to $\ddc[i^*(\pi_*(\Phi))]\wedge S$ in $\mD^0_{k-p+1}(X)$,  as $\lambda$ tends to infinity.
\end{lemma}
\proof
By Lemma \ref{l:Psi-lambda-lim}, we have that $\ddc \Psi_\lambda$ converges to $\ddc[i^*(\pi_*(\Phi))]\wedge S$ weakly. So the second assertion in the lemma is a direct consequence of the first one. We prove now the first assertion. We only consider $|\lambda|\geq 1$. 
We have by definition of $\Psi_\lambda$ 
$$\ddc\Psi_\lambda=(\pi_1)_*\big[A_\lambda^*(\ddc\Phi)\wedge \pi_2^*(S)\big].$$

\smallskip\noindent
{\bf Claim. } There is a positive closed current $\Theta_\lambda$ of bi-degree $(2k-p,2k-p)$ on $X\times X$ such that $\|\Theta_\lambda\|$ is bounded by a constant independent of $\lambda$ and $A_\lambda^*(\ddc\Phi)\wedge \pi_2^*(S)\geq -\Theta_\lambda$. 

\medskip

We first assume the claim and complete the proof of the lemma. Write $\ddc\Psi_\lambda=\Theta_\lambda^+-\Theta_\lambda^-$ with $\Theta_\lambda^-:=(\pi_1)_*(\Theta_\lambda)$ and $\Theta_\lambda^+:=\ddc\Psi_\lambda + \Theta_\lambda^-$. Clearly, $\Theta_\lambda^-$ is positive closed with mass bounded independently of $\lambda$. Since $\ddc\Psi_\lambda$ is an exact current, $\Theta_\lambda^+$ is closed and cohomologous to $\Theta_\lambda^-$. By the claim, this current is also positive. Since the mass of a positive closed current depends only on its cohomology class, the mass of $\Theta_\lambda^+$ is equal to the one of $\Theta_\lambda^-$ and therefore, is bounded independently of $\lambda$. The lemma follows.

We prove now the claim. Observe that $\Omega:=\ddc\|z'\|^2+\ddc\|w'\|^2$ is a strictly positive $(1,1)$-form on $W$.
Multiplying $\Phi$ by a constant, we can assume that
$\ddc\Phi\geq -\Omega^{2k-p-q}$ on $W$. It follows that 
$$A_\lambda^*(\ddc\Phi) \geq -(A_\lambda^*(\Omega))^{2k-p-q}=-(\ddc\|z'\|^2+|\lambda|^2\ddc \|w'\|^2)^{2k-p-q} \quad \text{on}\quad A_\lambda^{-1}(W).$$
Arguing as in \cite[Lemma 3.1]{DinhSibony15}, there is a sequence of smooth positive $(1,1)$-forms $\Omega_0,\Omega_1,\ldots$  on $X\times X$ such that 
$\sum \|\Omega_n\|$ is bounded by a constant independent of $\lambda$ and 
$$\sum \Omega_n\geq \ddc\|z'\|^2+|\lambda|^2\ddc \|w'\|^2 \quad \text{on} \quad  A_\lambda^{-1}(W)$$
(we can choose here $\Omega_0$ larger than $\ddc\|z'\|^2$). 

Recall again that the mass of a positive closed current depends only on its cohomology class. Thus, the current 
$\Theta_\lambda:= (\sum\Omega_n)^{2k-p-q}\wedge \pi_2^*(S)$ is well-defined and its mass is bounded by a constant independent of $\lambda$. We clearly have  $A_\lambda^*(\ddc\Phi)\wedge \pi_2^*(S)\geq -\Theta_\lambda$. 
This ends the proof of the claim.
\endproof

\noindent
{\bf End of the proof of Theorem \ref{t:main_2}.}
We need to show that $\mR=\pi^*(i_*(T\wedge S))$, where $T\wedge S$ is defined in the introduction using a super-potential of $T$. 
Let $\Phi$ be a real smooth test form of bi-degree $(2k-p-q,2k-p-q)$ with compact support on $\E$. We need to show that 
\begin{equation} \label{e:t2}
\langle \mR,\Phi\rangle =\big\langle \pi^*(i_*(T\wedge S)),\Phi\big\rangle. 
\end{equation}
Using a partition of unity and the notations as above, we can assume that $\Phi$ has compact support in $U\times \C^k$. Moreover, since both $\mR$ and $\pi^*(i_*(T\wedge S))$ are invariant under the action of $A_\lambda$, it is enough to consider the case where the support of $\Phi$ is contained in $W$ as in the situation of the above lemmas.

Let $\mU_T, U_T$ and $\alpha$ be as in the introduction. By the hypothesis of the theorem, $\mU_T$ is continuous on $\mD_{k-p+1}^0(X)$. 
Using this property,  \eqref{e:Psi-lambda} and Lemma \ref{l:ddc-Psi}, we obtain 
\begin{eqnarray*}
\langle \mR,\Phi\rangle & =  & \lim_{j\to\infty} \langle T,\Psi_{\lambda_{n_j}}\rangle= \lim_{j\to\infty} \langle T-\alpha,\Psi_{\lambda_{n_j}}\rangle + \big\langle \alpha ,i^*(\pi_*(\Phi))\wedge S \big\rangle \\
& =  & \lim_{j\to\infty} \langle U_T,\ddc \Psi_{\lambda_{n_j}}\rangle + \big\langle \alpha\wedge S ,i^*(\pi_*(\Phi))\big\rangle \\
& =  & \lim_{j\to\infty} \mU_T(\ddc \Psi_{\lambda_{n_j}}) + \big\langle \alpha \wedge S ,i^*(\pi_*(\Phi)) \big\rangle \\
& =  & \mU_T\big(\ddc [i^*(\pi_*(\Phi))]\wedge S\big) + \big\langle \alpha\wedge S ,i^*(\pi_*(\Phi)) \big\rangle.
\end{eqnarray*}
By the definition of wedge-product of currents, the last sum is equal to $\big\langle T\wedge S, i^*(\pi_*(\Phi)) \big\rangle$ which is equal to the RHS of \eqref{e:t2}. This ends the proof of the theorem.
\hfill $\square$

%%%%%%%%%%%%%%%%%%%%%%%%%%%%%%%%%%
\section{Holomorphic correspondences and actions on currents}  \label{s:dynamics}

As mentioned in the introduction, we will prove Theorem \ref{t:main_3} in a more general setting of holomorphic correspondences, see Theorem \ref{t:main_4} below. In this section, we first recall
some basic notions on correspondences.  We then study the action of holomorphic correspondences on currents and cohomology. Several notations that we use below,
have already been  defined in the introduction.

Consider an effective $k$-cycle $\Gamma =\sum \Gamma_j$ which is a finite combination of  irreducible analytic sets $\Gamma_j$ of dimension $k$ in $X\times X$.  We only consider $\Gamma_j$ such that its projections by $\pi_1$ and $\pi_2$ are equal to $X$. 
Note that the $ \Gamma_j$'s are not necessarily distinct. 
We say that $\Gamma$ 
defines a {\it (dominant) meromorphic correspondence} $f$ from $X$ to $X$ and $\Gamma$ is the {\it graph} of $f$.
If $A$ is  any subset of $X$,  define
$$f(A):=\pi_2(\pi_1^{-1}(A))\quad\text{and}\quad  f^{-1}(A):=\pi_1(\pi_2^{-1}(A)).$$
The {\it adjoint correspondence} of $f$ is denoted by $f^{-1}$. This is the  correspondence  whose graph
 is symmetric to $\Gamma$ with respect to the diagonal $\Delta$ of $X\times X$, that is, the  graph of $f^{-1}$ is  the image of $\Gamma$ by the involution $(x,y) \mapsto (y, x)$.

Define the {\it indeterminacy set} of $f$ by 
$$ I(f ) := \big\{x \in X, \quad \dim \pi_1^{-1} (x)\cap \Gamma > 0\big\}.$$
This is an analytic subset of co-dimension at least 2 of $X$. 
When $I(f)$ is empty, we say that $f$ is a {\it holomorphic correspondence}. 
Note that when $\pi_1$ restricted to $\Gamma$ is generically one-to-one, then
$f$ is a {\it meromorphic map} and when it is one-to-one, $f$ is a {\it holomorphic map}.
So we can consider correspondences as  multi-valued maps.  It is not difficult to show that if $f$ is a holomorphic map then $f^{-1}$ is a holomorphic correspondence but this is not true in general when $f$ is multi-valued, see Lemma \ref{l:inverse-map} below.

From now on, we only consider holomorphic correspondences, i.e. we always suppose $I(f)=\varnothing$. So for $x\in X$ the number of points in $f(x)$, counted with multiplicity, is a positive integer independent of $x$. We will denote it by $d_0$ or $d_0(f)$. Similarly, if $x$ is outside the indeterminacy set $I(f^{-1})$ of $f^{-1}$, then the number of points in $f^{-1}(x)$, counted with multiplicity, is a positive integer independent of $x$. We call it {\it topological degree} of $f$ and denote it by $d_k$ or $d_k(f)$, where $k$ is the dimension of $X$.
The degrees $d_0(f)$ and $d_k(f)$ coincide with the dynamical degrees that will be introduced later. If $g$ is another holomorphic correspondence on $X$, define $g\circ f$ as the correspondence such that $g\circ f(x)=g(f(x))$, counting multiplicity, for every $x\in X$. This is also a holomorphic correspondence. We have $d_0(g\circ f)=d_0(f)d_0(g)$ and $d_k(g\circ f)=d_k(d)d_k(g)$. As for holomorphic maps, the {\it iterate of order $n$} of $f$ is defined by $f^n:=f\circ\cdots \circ f$ ($n$ times). 
 
A holomorphic correspondence $f$ on $X$ induces pull-back and push-forward operators $f^*$ and $f_*$ on currents as follows : if $T$ is a current on $X$, define
$$f^* (T) := (\pi_1 )_* (\pi_2^* (T) \wedge [\Gamma]) \quad \text{and} \quad 
f_* (T) := (\pi_2 )_* (\pi_1^* (T) \wedge [\Gamma]),$$
when the wedge-products in these expressions are meaningful.
Clearly, when $T$ is smooth, the above currents $f^*(T)$ and $f_*(T)$ are well-defined.
These operators are in fact defined in more general settings. The following identity holds at least when $T$ and $\phi$ are smooth currents of the right degrees
\begin{equation} \label{e:pullback}
\langle f^*(T),\phi\rangle =\langle T,f_*(\phi)\rangle.
\end{equation}

Observe that the operators $f^*$ and $f_*$ commute with the operators $\partial$, $\dbar$ and recall that the Hodge cohomology groups can be defined using either smooth forms or singular currents. Therefore, $f^*$ and $f_*$ induce linear self-maps on the cohomology groups $H^{p,q}(X,\C)$ and $H^{p,p}(X,\R)$ that we will denote by the same notations. We deduce from \eqref{e:pullback} that 
$$f^*(c)\smallsmile c' = c\smallsmile f_*(c')$$
for all $c\in H^{p,q}(X,\C)$ and  $c'\in H^{k-p,k-q}(X,\C)$. Here, we identify the top degree group $H^{k,k}(X,\C)$ with $\C$ in the canonical way. So the last cup-products are complex numbers.
Thus, the operator $f^*$ on $H^{p,q}(X,\C)$ is dual to the operator $f_*$ on $H^{k-p,k-q}(X,\C)$.

\begin{proposition} \label{p:push-current}
Let $f$ be a holomorphic correspondence on $X$ as above and let $0\leq p\leq k$ be an integer. Then the operator $f_*$, which is well-defined on $\widetilde\mD_p(X)$, extends continuously to a linear operator from $\mD_p(X)$ to $\mD_p(X)$. Moreover, it preserves the subspace $\mD_p^0(X)$ and the cone $\mD_p^+(X)$ of positive $(p,p)$-currents. If $T$ is in $\mD_p(X)$, then 
$f_*\{T\}=\{f_*(T)\}$. 
If $g$ is another holomorphic correspondence on $X$, then $(g\circ f)_*=g_*\circ f_*$ on $\mD_p(X)$  and $H^{p,p}(X,\R)$. 
\end{proposition}
\proof
The proposition, except for the last assertion, is a particular case of \cite[Th.4.5]{DinhSibony07}. Also, by this theorem, if $T$ is in $\mD_p^+(X)$ and have no mass on proper analytic subsets of $X$, then $f_*(T)$ satisfies the same property. We will only prove the identity in the last assertion of the proposition for $\mD_p(X)$ because it implies the same identity for $H^{p,p}(X,\R)$. By continuity and Theorem \ref{t:reg}, we only need to check that $(g\circ f)_*(T)=g_*(f_*(T))$ for a smooth positive closed $(p,p)$-form $T$. Since the both sides of this identity are currents having no mass on proper analytic sets of $X$, it is enough to check the identity on a dense Zariski open subset of $X$. Note also that this identity is clear when $f$ and $g$ are maps, i.e. univalued. 

Let $\Gamma$ and $\Gamma'$ denote the graphs of $f$ and $g$ respectively. 
Choose a dense Zariski open set $\Omega$ of $X$ such that the restriction of $\pi_2$ (resp. $\pi_1$) to $\Gamma'\cap\pi_2^{-1}(\Omega)$ is a unramified covering (resp. unramified map, or equivalently, map without critical points). 
Define $\Omega':=g^{-1}(\Omega)$. Reducing $\Omega$ if necessary allows us to assume that 
the restriction of $\pi_2$ (resp. $\pi_1)$ to $\Gamma\cap\pi_2^{-1}(\Omega')$ is a unramified covering (resp. unramified map).  

Fix an arbitrary point $a$ in $\Omega$ and a small enough neighbourhood $U$ of $a$ in $\Omega$. Then the graph of $g$ restricted to $\pi_2^{-1}(U)$ is the union of the graphs of a finite number of holomorphic bijective maps, denoted by $g_i:U_i\to U$, which are defined on some open subsets $U_i$  of $\Omega'$. Since $U$ is small, $U_i$ is also small. Therefore, the graph of $f$ restricted to $\pi_2^{-1}(U_i)$ is the union of the graphs of a finite number of holomorphic bijective maps, denoted by $f_{ij}:U_{ij}\to U_i$, which are defined on some open subsets $U_{ij}$ of $\Omega''$. We deduce that the graph of $g\circ f$ restricted to  $\pi_2^{-1}(U)$ is the union of the graphs of $g_i\circ f_{ij}$. It follows that both $(g\circ f)_*(T)$ and $g_*(f_*(T))$ are equal in $U$ to the sum of $(g_i)_*(f_{ij})_*(T)$. This ends the proof of the proposition. 
\endproof

\begin{lemma} \label{l:bounded-sp}
Let $T$ be a positive closed $(p,p)$-current on $X$ with a bounded super-potential. Then $T$ has no mass on pluripolar sets of $X$. In particular, it has no mass on proper analytic subsets of $X$.
\end{lemma}
\proof
Let $E$ be a pluripolar set in $X$. Then there is a quasi-p.s.h. function $u$ on $X$ such that $u=-\infty$ on $E$. Subtracting from $u$ a constant allows us to assume that $u$ is strictly negative. A regularization theorem by Demailly says that we can find a sequence of negative smooth functions $u_n$ on $X$ decreasing to $u$ such that $\ddc u_n\geq -c\omega$ for some constant $c$ independent of $n$, see \cite{Demailly12}. We deduce that $\|\ddc u_n\|_*$ is bounded by a constant independent of $n$.

Let $\alpha$ and $\mU_T$ be as in the introduction. We have
$$\mU_T(\ddc u_n\wedge\omega^{k-p})= \langle T\wedge\omega^{k-p}, u_n\rangle - \langle \alpha \wedge\omega^{k-p}, u_n\rangle.$$
Therefore,
$$\lim_{n\to\infty}\mU_T(\ddc u_n\wedge\omega^{k-p})= \int_X u T\wedge\omega^{k-p} - \int_X u \alpha \wedge\omega^{k-p}.$$
By hypothesis, the LHS of the last identity is finite. Since quasi-p.s.h. functions are integrable, the second term in the RHS is finite. We deduce that the first term is also finite. So $u$ is integrable with respect to the trace measure of $T$. Since $u=-\infty$ on $E$, we deduce that $T$ has no mass on $E$.
\endproof

\begin{lemma} \label{l:pull-form}
Let $f$ and $p$ be as in Proposition \ref{p:push-current}.
Let $T$ be a smooth positive closed $(p,p)$-form. Then $f^*(T)$ is a positive closed $(p,p)$-current with a continuous super-potential.
\end{lemma}
\proof
For $R$ in $\mD_{k-p+1}^0(X)$ let $U_R$ be the potential of $R$ constructed in Proposition \ref{p:kernel}. Let $\alpha'$ be a smooth closed $(p,p)$-form on $X$ which is cohomologous to $f^*(T)$. Let $\mU$ denote the super-potential of $f^*(T)$ associated to the reference form $\alpha'$. When $R$ is smooth, we have $\mU(R)=\langle f^*(T),U_R\rangle -\langle \alpha',U_R\rangle$. Since $U_R$ depends continuously on $R$, the pairing $\langle \alpha',U_R\rangle$ extends to a continuous function of $R\in \mD_{k-p+1}^0(X)$. 
So, in order to get the lemma, it is enough to show that the pairing $\langle f^*(T),U_R\rangle$ satisfies the same property. 

We have $\langle f^*(T),U_R\rangle=\langle T, f_*(U_R)\rangle$ when $R$ is smooth. Since $T$ is smooth, it is enough to show that $f_*(U_R)$ is well-defined for every $R$ in $\mD_{k-p+1}^0(X)$, not necessarily smooth, and this current depends continuously on $R$. But this is a direct consequence of \cite[Th.4.6]{DinhSibony07}. The lemma follows.
\endproof

Let $\mD_p^{+c}(X)$ denote the cone of positive closed $(p,p)$-currents with continuous super-potentials. 
Let $\mD_p^c(X)$ be the vector subspace of $\mD_p(X)$ spanned by  $\mD_p^{+c}(X)$. It contains $\widetilde \mD_p(X)$. 
Define also $\mD_p^{0c}(X):=\mD_p^c(X)\cap \mD_p^0(X)$.
Observe that the notion of super-potential introduced in the introduction can be extended by linearity to all currents in $\mD_p(X)$. 
In particular, Lemma \ref{l:pull-form} holds for all (smooth) $T$ in $\widetilde\mD_p(X)$. 
If $T$ is a current in $\mD_p^0(X)$, we will use the super-potential of $T$ associated to the reference form $\alpha=0$; we call it {\it canonical} super-potential.

\begin{proposition} \label{p:pull-current}
Let $f$ and $p$ be as in Proposition \ref{p:push-current}. Then the operator $f^*$, which is well-defined on $\widetilde\mD_p(X)$, extends to a linear operator from $\mD_p^c(X)$ to itself with the following properties: 
\begin{enumerate}
\item[(i)]  The operator $f^*$ preserves the subspace $\mD_p^{0c}(X)$, the cone $\mD_p^{+c}(X)$ and we have $\{f^*(T)\}= f^*\{T\}$ for every $T$ in $\mD_p^c(X)$;
\item[(ii)]  Let $T$ be a current in $\mD_p^{0c}(X)$. If $\mU_T$ and $\mU_{f^*(T)}$ denote the canonical super-potentials of $T$ and $f^*(T)$ respectively, then $\mU_{f^*(T)}(R)=\mU_T(f_*(R))$ for every $R\in\mD^0_{k-p+1}(X)$.  
\end{enumerate}
\end{proposition}
\proof
Let $\Gamma$ denote the graph of $f$ in $X\times X$ and $[\Gamma]$ the positive closed $(k,k)$-current of integration on $\Gamma$.
Let $T$ be any current in $\mD_p^{+c}(X)$. By Lemma \ref{l:sp_pullback}, $\pi_2^*(T)$ is a current in $\mD_p^{+c}(X\times X)$ and then by Theorem \ref{t:main_2}, the wedge-product $\pi_2^*(T)\wedge S$ defines a positive closed current for every positive closed $(k,k)$-current $S$ on $X\times X$, see also \cite{Vu}. This wedge-product depends continuously on $S$ and therefore, we easily deduce from Theorem \ref{t:reg} that 
$$\{\pi_2^*(T)\wedge S\}=\{\pi_2^*(T)\}\smallsmile \{S\}=\pi_2^*\{T\}\smallsmile\{S\}$$
because this obviously holds when $S$ is smooth.
In particular, we can define
$$f^*(T):=(\pi_1)_*(\pi_2^*(T)\wedge[\Gamma])$$
which is  a positive closed $(p,p)$-current. If $T$ is smooth, this is the usual pull-back of $T$ by $f$ mentioned in the beginning of the section. We also have 
$$\{f^*(T)\} = (\pi_1)_*\big(\pi_2^*\{T\}\smallsmile\{\Gamma\}\big)=f^*\{T\}.$$
By linearity, $f^*$ extends to $\mD_p^c(X)$ and the last identity holds for the extended operator.

In order to prove that $f^*$ preserves $\mD_p^c(X),$ $ \mD^{0c}_p(X)$ and $\mD^{+c}_p(X),$  
we claim now that it is enough to prove the identity in (ii) for $R$ smooth. Indeed, by Proposition \ref{p:push-current}, the RHS of this identity extends to
a continuous function on $R\in\mD_{k-p+1}^0(X)$. So the identity implies that $f^*(T)$ has a continuous super-potential and therefore, Lemma \ref{l:pull-form} implies the 
same property for every current $T$ in $\mD_p^c(X)$ because we can write it as the sum of a current in $\mD_p^{0c}(X)$ and a smooth closed form.  This proves (i).

We prove now the identity in (ii) for $R$ smooth. 
Let $U_R$ be as in the proof of Lemma \ref{l:pull-form}.
We have 
$$\mU_{f^*(T)}(R)=\langle f_*(T), U_R\rangle =\langle \pi_2^*(T)\wedge [\Gamma], \pi_1^*(U_R)\rangle.$$
By Theorem \ref{t:reg}, there is a sequence of smooth forms $[\Gamma]_n$ converging to $[\Gamma]$ in $\mD_k(X\times X)$. 
We deduce from the last identities and the continuity of the wedge-product that 
$$\mU_{f^*(T)}(R)=\lim_{n\to\infty} \big\langle \pi_2^*(T)\wedge [\Gamma]_n, \pi_1^*(U_R)\big\rangle=\lim_{n\to\infty} \big\langle T, (\pi_2)_*([\Gamma]_n\wedge \pi_1^*(U_R))\big\rangle.$$
The last pairing is equal to the value of $\mU_T$ at $\ddc (\pi_2)_*([\Gamma]_n\wedge \pi_1^*(U_R))$. Observe that the last current is equal to $(\pi_2)_*([\Gamma]_n\wedge \pi_1^*(R))$ and therefore, it 
 converges to $f_*(R)$ in $\mD^0_{k-p+1}(X)$.
Since  $\mU_T$ is continuous, we deduce that the last limit is equal to $\mU_T(f_*(R))$. This ends the proof of the proposition.
\endproof

\begin{remark} \rm \label{r:PC-PB}
Let $T, R$ be as in Proposition \ref{p:pull-current}(ii) and let $\mU_R$ denote the canonical super-potential of $R$. When $T$ is smooth, we have $\mU_R(T)=\langle U_R,T\rangle=\mU_T(R)$. Therefore, we can extend $\mU_R$ to $\mD_p^{0c}(X)$ by setting $\mU_R(T):=\mU_T(R)$. We deduce from Proposition \ref{p:pull-current}(ii) that $\mU_{f_*(R)}(T)=\mU_R(f^*(T))$ for $T\in \mD_p^{0c}(X)$.
\end{remark}

Let $\Omega$ be a dense Zariski open set in $X$ so that the restriction $\tau_1$ (resp. $\tau_2$) of $\pi_1$ (resp. $\pi_2$) to $\Gamma\cap\pi_1^{-1}(\Omega)$  is a unramified covering (resp. unramified map). So $\tau_1$ and $\tau_2$ are locally bi-holomorphic maps. 
If $T$ is any current on $X$, we can define the current $(\tau_1)_*\tau_2^*(T)$ on $\Omega$ which depends continuously on $T$.
The following lemma gives us a more explicit description of the operator $f^*$.

\begin{lemma} \label{l:pull-explicit}
Let $f,g,p$ be as in Proposition \ref{p:push-current} and 
 $\Gamma, \Omega,\tau_1,\tau_2$ be as above. Let $T$ be a current in $\mD_p^c(X)$. Then $f^*(T)$ is the extension by $0$ of 
 $(\tau_1)_*\tau_2^*(T)$ to $X$. Moreover, the operator $f^*$ satisfies the following continuity property : 
 if $(T_n)_{n\geq 0}$ is a sequence in   $\mD_p^{+c}(X)$, converging weakly to a current $T$ also in $\mD_p^{+c}(X)$, then $f^*(T_n)$ converges weakly to $f^*(T)$. 
We also have $(g\circ f)^*=f^*\circ g^*$ on $\mD_p^c(X)$.
\end{lemma}
\proof
By Proposition \ref{p:pull-current}, $f^*(T)$ has a continuous super-potential. 
By Lemma \ref{l:bounded-sp}, this current has no mass on $X\setminus\Omega$. So for the first assertion in the lemma, we only need to show that $f^*(T)=(\tau_1)_*\tau_2^*(T)$ on $\Omega$. Let $\phi$ be a real smooth $(k-p,k-p)$-form with compact support on $\Omega$. We need to check that $\langle f^*(T),\phi\rangle =\langle (\tau_1)_*\tau_2^*(T),\phi\rangle$. 
This identity is clear if we replace $T$ by a smooth form. Therefore, adding to $T$ a suitable smooth closed form, we can assume that $T$ belongs to $\mD_p^{0c}(X)$. By Proposition \ref{p:pull-current}, we have $\mU_{f^*(T)}(\ddc\phi)=\mU_T(\ddc f_*(\phi))$. 
Using the properties of $\Omega$, we see that $f_*(\phi)$ is smooth and is equal to $(\tau_2)_*\tau_1^*(\phi)$. Therefore, we have
$$\langle f^*(T), \phi\rangle =\mU_{f^*(T)}(\ddc \phi)  =\mU_T(\ddc f_*(\phi)) =\langle T, f_*(\phi)\rangle =\langle T, (\tau_2)_*\tau_1^*(\phi)\rangle.$$
This implies the desired identity and then the first assertion in the lemma.

We prove now the second assertion. Since $T_n$ converges to $T$, its mass and cohomology class are bounded independently of $n$. Therefore, the same properties hold for $f^*(T_n)$. Extracting a subsequence allows us to assume that $f^*(T_n)$ converges to a positive closed current $S$. We need to show that $S=f^*(T)$. Since the operator $(\tau_1)_*(\tau_2)^*$ is continuous, we deduce that $S=f^*(T)$ on $\Omega$.  Moreover, since $f^*(T)$ has no mass on $X\setminus\Omega$, we obtain that $S-f^*(T)$ is a positive closed current with support in $X\setminus \Omega$, see \cite{Skoda}. By Proposition \ref{p:pull-current}, the operator $f^*$ is compatible with the action of $f$ on cohomology. Thus, the cohomology class of $S-f^*(T)$ is 0 and therefore this positive closed current is 0. The result follows. 

Consider now the last assertion in the lemma. Let $\Gamma'$ be the graph of $g$ and define $\Omega':=f(\Omega)$. We can choose $\Omega$ so that $\pi_1$ restricted to $\Gamma'\cap \pi_1^{-1}(\Omega')$ is a unramified covering over $\Omega'$ and that $\pi_2$ restricted to $\Gamma'\cap \pi_1^{-1}(\Omega')$ is a unramified map. Using these maps and $\tau_1,\tau_2$, we can see as in the proof of Proposition \ref{p:push-current} that $(g\circ f)^*=f^*\circ g^*$ on $\Omega$. 
We then deduce from the first assertion of the lemma that $(g\circ f)^*=f^*\circ g^*$ on $\mD_p^c(X)$.
\endproof

By Proposition \ref{p:push-current}, the operator $(f^n)_*$ acts continuously on $\mD_p(X)$ and we have $(f^n)_*=(f_*)^n$. It follows that $(f^n)_*=(f_*)^n$ on the cohomology group $H^{p,p}(X,\R)$. Similarly, we deduce from Proposition \ref{p:pull-current} and Lemma \ref{l:pull-explicit} that  $(f^n)^*=(f^*)^n$ on $\mD_p^c(X)$ and on $H^{p,p}(X,\R)$. We call {\it dynamical degree} of order $p$ of $f$ the spectral radius $d_p(f)$ of $f^*$ on $H^{p,p}(X,\R)$. Since the last operator is dual to $f_*$ acting on $H^{k-p,k-p}(X,\R)$, $d_p(f)$ is also the spectral radius of $f_*$ on $H^{k-p,k-p}(X,\R)$. It is not difficult to show that this dynamical degree is an eigenvalue of $f^*$ on $H^{p,p}(X,\R)$ and of  $f_*$ on $H^{k-p,k-p}(X,\R)$. We have the following elementary lemma.

\begin{lemma} \label{l:inverse-map}
Let $f$ be a surjective holomorphic map from $X$ to $X$. Then $f^{-1}$ is a holomorphic correspondence.
\end{lemma}
\proof
Let $a$ be a point in $X$. We need to show that $f^{-1}(a)$ is finite. Assume by contradiction that $f^{-1}(a)$ contains an analytic subset $Y$ of dimension $q\geq 1$. We can assume that $Y$ is of pure dimension $q$. Since the positive closed current  $(k-q,k-q)$-current defined by $Y$ is not zero, its class in $H^{k-q,k-q}(X,\R)$ is not zero.  We denote this class by $\{Y\}$. The choice of $Y$ implies that the positive closed current $f_*([Y])$ is supported by the point $a$. We used here that $f$ is univalued. Therefore, it vanishes and we obtain $f_*\{Y\}=0$. 

Observe now that $f_*\circ f^*(\alpha)=d_k(f) \alpha$ for every smooth form $\alpha$. We deduce that $f_*\circ f^*$ acting on $H^{k-q,k-q}(X,\R)$ is just the multiplication by $d_k(f)$. In particular, the operator $f_*$ is invertible. But $f_*\{Y\}=0$. This is a contradiction. The lemma follows.
\endproof

%%%%%%%%%%%%%%%%%%%%%%%%%%%%%%%%%%
\section{Dynamical Green currents and periodic points} \label{s:per} 

In this section, we will prove Theorem \ref{t:main_4} below which is  an extension of Theorem \ref{t:main_3} for correspondences. Throughout the section, we assume that $f$ is a correspondence as in this statement.

\begin{theorem} \label{t:main_4}
Let $f$ be a  (dominant) holomorphic correspondence on a compact K\"ahler manifold $X$ of dimension $k$ such that $f^{-1}$ is also a holomorphic correspondence. Assume that the action of $f$ on cohomology is simple. More precisely, there is an integer $0\leq p\leq k$ such that (1) the dynamical degree of order $p$ of $f$, denoted by $d$,  is strictly larger than the other dynamical degrees, (2) $d$ is a simple eigenvalue  of $f^*$ on $H^{p,p}(X,\R)$,
and (3) the other (real or complex) eigenvalues of this operator have modulus strictly smaller than $d$. If $P_n$ denotes the number of isolated periodic points of period $n$ of $f$, counted with multiplicities, then we have
$P_n\leq d^n+o(d^n)$ as $n$ goes to infinity.
\end{theorem}

We call $d$ the {\it main dynamical degree} of $f$.  Note that when $\dim H^{p,p}(X,\R)=1$, e.g. when $p=0$ or $p=k$, the properties (2) and (3) are automatically satisfied. 
The quantity $h_a(f):=\log d$ is the {\it algebraic entropy} of the correspondence.  
Let $\Gamma_n$ denote the graph of $f^n$ in $X\times X$. A point $x$ is an {\it isolated periodic point} of period $n$ of $f$ if there is at least a germ of analytic set of
$\Gamma_n$ at $(x,x)$ which intersects $\Delta$ only at this point $(x,x)$. The {\it multiplicity} of this periodic point is the total multiplicity of the intersections of such germs of analytic sets with $\Delta$ at $(x,x)$.   

In order to prove Theorem \ref{t:main_4}, we will use the property that the sequence $\Gamma_n$, suitably normalized, admits a limit in the sense of currents when $n$ tends to infinity, see Proposition \ref{p:graph} below. We first need to construct and study the two main dynamical Green currents $T^+$ and $T^-$ of $f$, in particular, the action of $f^n$ on cohomology. 

\begin{lemma} \label{l:Green+class}
The sequence $d^{-n} (f^n)^*\{\omega^p\}$ converges  in $H^{p,p}(X,\R)$ to a non-zero class $c^+$. Let $L^+$ denote the real line spanned by $c^+$. Then we can write $H^{p,p}(X,\R)=L^+\oplus H^+$, where $H^+$ is a hyperplane of $H^{p,p}(X,\R)$ which is invariant by $f^*$. Moreover, we have $d^{-1}f^*(c^+)=c^+$ and the spectral radius of $d^{-1}f^*$ on $H^+$ is strictly smaller than $1$. 
\end{lemma}
\proof
Observe that 1 is a simple eigenvalue of $d^{-1}f^*$ on $H^{p,p}(X,\R)$ and is the only eigenvalue of maximal modulus 1. 
Let $L^+$ denote the set of vectors which are fixed by $d^{-1} f^*$. Clearly, $L^+$ is a real line. It is easy to obtain from basic results of linear algebra 
that there is a hyperplane $H^+$, corresponding to the other eigenvalues of $d^{-1}f^*$, which is invariant by this operator. We have  $H^{p,p}(X,\R)=L^+\oplus H^+$ and the spectral radius of $d^{-1}f^*$ on $H^+$ is strictly smaller than $1$.  

Write $\{\omega^p\}=c^++c_0$, where $c^+$ is a class in $L^+$ and $c_0$ is a class in $H^+$. We have  $d^{-n}(f^n)^*\{\omega^p\}=c^++d^{-n}(f^n)^*(c_0)$ and therefore $d^{-n}(f^n)^*\{\omega^p\}$ converges to $c^+$ since on $H^+$ the sequence $d^{-n}(f^n)^*$ converges to 0. It remains to check that $c^+\not=0$. 
We will prove a more general property in Lemma \ref{l:c+c-} below.
\endproof

\begin{lemma} \label{l:Green-class}
The sequence $d^{-n} (f^n)_*\{\omega^{k-p}\}$ converges in $H^{k-p,k-p}(X,\R)$  to a non-zero class $c^-$. Let $L^-$ be the real line spanned by $c^-$. Then we can write $H^{k-p,k-p}(X,\R)=L^-\oplus H^-$, where $H^-$ is a hyperplane of $H^{k-p,k-p}(X,\R)$, invariant by $f_*$. Moreover, we have $d^{-1}f_*(c^-)=c^-$ and the spectral radius of $d^{-1}f_*$ on $H_-$ is strictly smaller than $1$. 
\end{lemma}
\proof
The proof is exactly the same as for Lemma \ref{l:Green+class}. We also obtain that $c^-\not=0$ as a consequence of Lemma \ref{l:c+c-} below. 
\endproof

\begin{lemma} \label{l:c+c-}
The intersection number  $c^+\smallsmile c^-$ is strictly positive. In order to simplify the notation, we can (and we do) multiply $\omega$ by a constant so that 
$$c^+\smallsmile c^-=c^+\smallsmile \{\omega^{k-p}\}=\{\omega^p\}\smallsmile c^-=1.$$ 
\end{lemma}
\proof
Using the definitions of $c^+$, $c^-$ and that the operator $f^*$ on $H^{p,p}(X,\R)$ is dual to $f_*$ on $H^{k-p,k-p}(X,\R)$, we obtain
$$c^+\smallsmile c^-=\lim_{n\to\infty} d^{-2n} (f^n)^*\{\omega^p\}\smallsmile (f^n)_*\{\omega^{k-p}\}=\lim_{n\to\infty} d^{-2n} (f^{2n})^*\{\omega^p\}\smallsmile \{\omega^{k-p}\}.$$
It follows that $c^+\smallsmile c^-=c^+\smallsmile \{\omega^{k-p}\}$ and similarly, we get $c^+\smallsmile c^-=\{\omega^p\}\smallsmile c^-$. So if these cup-products are positive numbers, we can multiply $\omega$ by a constant in order to assume that they are all equal to 1.

The cup-product in the last limit is equal to the mass of the positive closed current $d^{-2n} (f^n)^*(\omega^p)$. Observe that 
any limit value $T$ of the sequence  $d^{-2n} (f^{2n})^*(\omega^p)$ is a positive closed current in the class $c^+$
and the mass of $T$ is equal to $c^+\smallsmile c^-$. Therefore, this is a non-negative number and we only need to check that $T$ is not equal to 0.

Observe that if $c$ is a class in $H^{p,p}(X,\C)\setminus H^+$ then $d^{-n}(f^n)^*(c)$ converges to a non-zero class in $L^+$. 
Since any class $c$ can be written as the difference of two classes of positive smooth forms, there is a positive closed  $(p,p)$-smooth form $\beta$ such that $\{\beta\}$ is outside $H^+$ and therefore $d^{-n}(f^n)^*\{\beta\}$ converges to a class $c'$ in $L^+\setminus\{0\}$.
Multiplying $\beta$ with a positive constant allows us to assume that $\beta\leq \omega^p$. It follows that there is a limit value of 
$d^{-2n} (f^{2n})^*(\beta)$ which is a positive closed current $T'$ smaller than or equal to $T$. Since $T'$ belongs to the class $c'$, it is not equal to 0. It follows that $T$ is not equal to 0. This ends the proof of the lemma.
\endproof

\begin{proposition} \label{p:Green}
The sequence $d^{-n}(f^n)^*(\omega^p)$ converges weakly to a positive closed $(p,p)$-current $T^+$ in the cohomology class $c^+$, as $n$ tends to infinity. The sequence $d^{-n}(f^n)_*(\omega^{k-p})$ converges weakly to a positive closed $(k-p,k-p)$-current $T^-$ in the cohomology class $c^-$, as $n$ tends to infinity. Moreover, $T^+$ and $T^-$ have continuous super-potentials.  
\end{proposition}

We will only prove the first assertion. The second one is obtained in the same way by using $f^{-1}$ instead of $f$. 
Consider the operator $\Lambda:=d^{-1}f_*$ acting on  $\mD^0_{k-p+1}(X)$. 
By Proposition \ref{p:push-current}, $\Lambda$ is bounded and continuous with respect to the topology considered on $\mD^0_{k-p+1}(X)$.
Recall that the action of $f^*$ on cohomology of $X$ is simple and that the mass of a positive closed current depends only on its cohomology class. Therefore, there is a constant $0<\delta<1$ such that the norm of $\Lambda^n$ with respect to the $\ast$-norm is bounded by $O(\delta^n)$.

Choose a smooth real closed $(p,p)$-form $\alpha_0$ in the class $c^+$. 
Define  $\alpha_n:=d^{-n}(f^n)^*(\alpha_0)$. This current is in the cohomology class $c^+$. In particular, 
the class of $\alpha_1-\alpha_0$ is zero and we denote by $\mV$ its canonical super-potential. 
By Proposition \ref{p:pull-current}, this super-potential is continuous on $\mD_{k-p+1}^0(X)$. More precisely, $\mV$ is bounded on $\ast$-bounded subsets of $\mD_{k-p+1}^0(X)$ and is continuous with respect to the topology considered on this space.

\begin{lemma} \label{l:pull-alpha}
The canonical super-potential of $\alpha_n-\alpha_0$, denoted by $\mV_n$, is continuous on $\mD_{k-p+1}^0(X)$ and is equal to
$\mV+\mV\circ\Lambda +\cdots+\mV\circ\Lambda^{n-1}$. Moreover, it converges to a continuous function $\mV_\infty$ on $\mD_{k-p+1}^0(X)$ as $n$ goes to infinity.
\end{lemma}
\proof
We have 
$$\alpha_n-\alpha_{n-1} = d^{-n+1}(f^{n-1})^*(\alpha_1-\alpha_0).$$
By Proposition \ref{p:pull-current},  the canonical super-potential of $\alpha_n-\alpha_{n-1}$ is equal to $\mV\circ \Lambda^{n-1}$. By taking a sum of such currents, we obtain that 
$$\mV_n=\mV+\mV\circ\Lambda +\cdots+\mV\circ\Lambda^{n-1}.$$
Since $\Lambda$ is bounded and continuous, $\mV_n$ is a continuous super-potential.

Since the $\ast$-norm of $\Lambda^n$ is bounded by $O(\delta^n)$. The lass sum converges, uniformly on $\ast$-bounded sets in $\mD_{k-p+1}^0(X)$,  to a function $\mV_\infty$ on $\mD_{k-p+1}^0(X)$. We have 
$$\mV_\infty=\mV+\mV\circ\Lambda +\cdots+\mV\circ\Lambda^n+\cdots$$
Clearly, this function is bounded on each $\ast$-bounded subset of $\mD_{k-p+1}^0(X)$.
On each $\ast$-bounded subset of $\mD_{k-p+1}^0(X)$, when $n$ is large enough, the partial sum of order $n$ of the last series is continuous and the rest is small. It is not difficult to deduce that $\mV_\infty$ is continuous  with respect to the topology considered on $\mD_{k-p+1}^0(X)$.
\endproof

Fix a real vector space $\Hc$ spanned by a finite number of smooth real closed $(p,p)$-forms such that the map $\alpha\mapsto \{\alpha\}$ defines a bijection from $\Hc$ to 
the hyperplane $H^+$ in $H^{p,p}(X,\R)$. Define $\Omega_0:=\omega^p-\alpha_0$. By definition of $c^+$, the class $\{\Omega_0\}$ is in $H^+$. Therefore, we can choose $\Hc$ so that $\Omega_0$ is in $\Hc$. 

For each class $c$ in $H^+$, there is a unique form $\beta_c$ in $\Hc$ such that $\{\beta_c\}=c$. If $c':=d^{-1}f^*(c)$, then the class of $d^{-1}f^*(\beta_c)-\beta_{c'}$ is 0 and we denote by $\mW_c$ its canonical super-potential. We see that the set of all $\mW_c$ is a finite dimensional vector space because $\mW_c$ depends linearly on $c$.
Let $c_0$ denote the class of $\Omega_0$ and define $c_n:=d^{-n}(f^n)^*(c_0)$. 
For simplicity, denote by $\beta_n$ the unique form in $\Hc$ such that $\{\beta_n\}=c_n$ and define $\mW_n:=\mW_{c_n}$. Define also $\Omega_n:=d^{-n}(f^n)^*(\Omega_0)$ and $\mU'_n$ the canonical super-potential of $\Omega_n-\beta_n$.
By Proposition \ref{p:pull-current}, all $\mW_c, \mW_n$ and $\mU'_n$ are continuous. Note that $\Omega_0=\beta_0$ and hence $\mU'_0=0$. Note also that since the spectral radius of $d^{-1}f^*$ on $H^+$ is smaller than 1, fixing a norm for $H^{p,p}(X,\R)$, we have $\|c_n\|=O(\delta^n)$  and $\beta_n=O(\delta^n)$ for some constant $0<\delta<1$. 

\begin{lemma} \label{l:pull-Omega}
The sequence $\mU'_n$ converges to $0$ on $\mD_{k-p+1}^0(X)$ when $n$ tends to infinity. 
\end{lemma}
\proof
We have 
$$\Omega_n-d^{-1}f^*(\beta_{n-1})=d^{-1} f^*(\Omega_{n-1}-\beta_{n-1}).$$
Therefore, the canonical super-potential of $\Omega_n-d^{-1}f^*(\beta_{n-1})$ is equal to $\mU'_{n-1}\circ\Lambda$.
We then deduce from the above discussion that $\mU'_n=\mU'_{n-1}\circ\Lambda + \mW_{n-1}$. By induction and using that $\mU'_0=0$, we obtain 
$$\mU'_n:=\mW_0\circ\Lambda^{n-1} +\mW_1\circ \Lambda^{n-2}+\cdots +\mW_{n-1}.$$

Since $\mW_c$ depends linearly on $c$, there is a constant $A>0$ such that for all $c\in H^+$ and $R\in \mD_{k-p+1}^0(X)$, we have $|\mW_c(R)|\leq A \|c\|\|R\|_*$. Recall that  $\|c_n\|=O(\delta^n)$ and $\|\Lambda^n\|=O(\delta^n)$ for some constant $0<\delta<1$. We then deduce that 
$$|\mW_m(\Lambda^{n-m-1}(R))|\leq A \|c_m\|\|\Lambda^{n-m-1}(R)\|_* \leq A'\delta^n\|R\|_*$$
for some constant $A'>0$. It is now clear that $\mU'_n(R)$ converges to 0 as $n$ goes to infinity. This ends the proof of the lemma.
\endproof

\noindent
{\bf End of the proof of Proposition \ref{p:Green}.} 
Observe that $d^{-n}(f^n)^*(\omega^p)$ is a positive closed current with bounded mass. So any limit of $d^{-n}(f^n)^*(\omega^p)$, when $n$ goes to infinity, is a positive closed current. Fix any sequence $n_j$ going to infinity such that 
$d^{-n_j}(f^{n_j})^*(\omega^p)$ converges to some positive closed current $T^+$. Clearly,  $T^+$ is cohomologous to $\alpha_0$.
 Denote by $\mU$ the canonical super-potential of $T^+-\alpha_0$. 
 
We have $d^{-n}(f^n)^*(\omega^p)=\alpha_n+\Omega_n$. So $d^{-n}(f^n)^*(\omega^p)$ is cohomologous
to $\alpha_0+\beta_n$. Let $\mU_n$ denote the canonical super-potential of $d^{-n}(f^n)^*(\omega^p)-(\alpha_0+\beta_n)$.
We have $\mU_n=\mV_n+\mU_n'.$ Hence, by Lemmas \ref{l:pull-alpha} and \ref{l:pull-Omega}, $\mU_n$ converges to $\mV_\infty$ which is continuous on $\mD_{k-p+1}^0(X)$. So, it remains to check that  the canonical super-potential $\mU$ of $T^+-\alpha_0$ is equal to $\mV_\infty$ on $\widetilde \mD_{k-p+1}^0(X)$. Indeed, this property implies that $T^+$ does not depend on the choice of the sequence $n_j$ as $\mV_\infty$ satisfies the same property, or equivalently,  $d^{-n}(f^n)^*(\omega^p)$ converges to $T^+$, which is a current of  continuous super-potentials.

Fix a smooth form $R$ in $\widetilde\mD_{k-p+1}^0(X)$ and let $U_R$ be a smooth potential of $R$. Denote by $\mU_R$ the canonical super-potential of $R$ which is
a continuous super-potential because $R$ is smooth. We  have 
$$\mU_n(R) = \big\langle  d^{-n}(f^n)^*(\omega^p)-(\alpha_0+\beta_n), U_R\big\rangle = \mU_R\big(d^{-n}(f^n)^*(\omega^p)-(\alpha_0+\beta_n)\big).$$
Then, using the continuity of $\mU_R$, we obtain
$$\mU(R)=\mU_R(T^+-\alpha_0)=\lim_{j\to\infty} \mU_R\big(d^{-n_j}(f^{n_j})^*(\omega^p)-(\alpha_0+\beta_{n_j})\big)=\lim_{j\to\infty} \mU_{n_j}(R)=\mV_\infty(R).$$
This ends the proof of the proposition.
\hfill $\square$

\medskip

We will need the following lemma.

\begin{lemma} \label{l:conv-rs}
Let $\phi$ be a smooth $(r,s)$-form on $X$ with $(r,s)\not=(p,p)$. 
Then $d^{-n}(f^n)^*(\phi)$ converges weakly to $0$ as $n$ tends to infinity. Let $\psi$ be a smooth  $(r,s)$-form on $X$ with $(r,s)\not=(k-p,k-p)$.  Then $d^{-n}(f^n)_*(\psi)$ converges weakly to $0$ as $n$ tends to infinity.
\end{lemma}
\proof
We will prove the first assertion. The second one can be obtained in the same way using $f^{-1}$ instead of $f$. Let $\phi'$ be a smooth $(k-r,k-s)$-form. We need to show that $\big\langle d^{-n} (f^n)^*(\phi),\phi'\big\rangle$ tends to 0 as $n$ tends to infinity. 
We distinguish two cases.

\medskip\noindent
{\bf Case 1.}   Assume that $r=s\not=p$. We can for simplicity assume that $\phi$ and $\phi'$ are positive forms because they can be written as linear combinations of such forms. Multiplying them by a constant allows us to assume that $\phi\leq \omega^r$ and $\phi'\leq \omega^{k-r}$. We have
$$\big\langle d^{-n} (f^n)^*(\phi),\phi' \big\rangle \leq \big\langle d^{-n} (f^n)^*(\omega^r),\omega^{k-r} \big\rangle \leq d^{-n}\|(f^n)^*(\omega^r)\|.$$
Recall that the mass of a positive closed current depends only on its cohomology class. Moreover, since $r\not=p$, the spectral radius of $f^*$ on $H^{r,r}(X,\R)$ is strictly smaller than $d$. It follows that the sequence of operators $d^{-n}(f^n)^*$ on $H^{r,r}(X,\R)$ converges to 0 as $n$ tends to infinity.  Thus, $d^{-n}(f^n)^*(\omega^r)$ tends to 0. The result follows.

\medskip\noindent
{\bf Case 2.} Assume now that $r\not=s$. So we have $(r,s)\not=(p,p)$. Assume for simplicity that $r>s$. The opposite case can be treated in the same way. 
We can assume that $\phi=\theta\wedge \Omega$ and $\phi'=\theta'\wedge \Omega'$, where $\theta, \theta'$ are smooth forms of bi-degrees $(r-s,0), (0,r-s)$, and $\Omega,\Omega'$ are smooth positive forms of bi-degrees $(s,s), (k-r,k-r)$ respectively. Indeed, $\phi$ and $\phi'$ can be written as finite sums of forms of the considered types. Define $\widetilde\Omega:=\theta\wedge \overline\theta\wedge\Omega$ and $\widetilde\Omega':=\theta'\wedge \overline\theta'\wedge\Omega'$.
The form $\widetilde\Omega$ is of bi-degree $(r,r)$ and the form $\widetilde\Omega'$ is of bi-degree $(k-s,k-s)$.  
Using Cauchy-Schwarz's inequality, we get
\begin{eqnarray*}
\big|\big\langle d^{-n} (f^n)^*(\phi),\phi' \big\rangle \big|^2 & = & \big|\big\langle d^{-n}[\Gamma_n],\pi_2^*(\theta\wedge\Omega)\wedge\pi_1^*(\theta'\wedge\Omega')\big\rangle\big|^2 \\
& \leq &  \big|\big\langle d^{-n}[\Gamma_n],\pi_2^*(\widetilde\Omega)\wedge\pi_1^*(\Omega')\big\rangle\big|  \big|\big\langle d^{-n}[\Gamma_n],\pi_2^*(\Omega)\wedge \pi_1^*(\widetilde\Omega')\big\rangle\big| \\
& \leq &  \big|\big\langle d^{-n}(f^n)^*(\widetilde\Omega),\Omega'\big\rangle\big|  \big|\big\langle d^{-n}(f^n)^*(\Omega), \widetilde\Omega' \big\rangle\big| .
\end{eqnarray*}

The last product can be bounded by a constant times
$$ \big|\big\langle d^{-n}(f^n)^*(\omega^r),\omega^{k-r}\big\rangle\big|  \big|\big\langle d^{-n}(f^n)^*(\omega^s), \omega^{k-s}  \big\rangle\big| = \|d^{-n}(f^n)^*(\omega^r)\| \|d^{-n}(f^n)^*(\omega^s)\|.$$
As in Case 1, we see that both factors in the last product are bounded. Moreover, since $(r,s)\not=(p,p)$, at least one of these factors converges to 0 as $n$ goes to infinity. 
We deduce that $\big\langle d^{-n} (f^n)^*(\phi),\phi' \big\rangle$ tends to 0. This ends the proof of 
the lemma.
\endproof

\begin{proposition} \label{p:Green-eq}
Let $\phi$ be a smooth $(p,p)$-form on $X$ and $\lambda_\phi:=\langle T^-,\phi\rangle$. 
Then $d^{-n}(f^n)^*(\phi)$ converges weakly to $\lambda_\phi T^+$ as $n$ tends to infinity. Let $\psi$ be a smooth  $(k-p,k-p)$-form on $X$ and $\lambda_\psi:=\langle T^+,\psi\rangle$. 
Then $d^{-n}(f^n)_*(\psi)$ converges weakly to $\lambda_\psi T^-$ as $n$ tends to infinity.
\end{proposition}
\proof
We only prove the first assertion. The second one can be obtained in the same way using $f^{-1}$ instead of $f$.
It is enough to consider the case where $\phi$ is positive because $\phi$ can be written as a linear combination of positive smooth forms. Multiplying $\phi$ by 
a constant allows us to assume that $\phi\leq \omega^p$. It follows that  $d^{-n}(f^n)^*(\phi)\leq d^{-n}(f^n)^*(\omega^p)$. 
By Proposition \ref{p:Green}, $d^{-n}(f^n)^*(\omega^p)$ converges to $T^+$. Therefore, the positive currents $d^{-n}(f^n)^*(\phi)$ have bounded masses and 
this family of currents is relatively compact, i.e., it admits convergent subsequences. 

Let $\Fc$ denote the family of the limits of all such convergent subsequences. This is a compact family of positive $(p,p)$-currents which are bounded by $T^+$. We also have $\partial d^{-n}(f^n)^*(\phi) = d^{-n}(f^n)^*(\partial \phi)$. Lemma \ref{l:conv-rs} implies that the last expression converges weakly to 0. So all currents in $\Fc$ are $\partial$-closed. Similarly, we obtain that they are $\dbar$-closed and hence closed. By Theorem \ref{t:main_1} and Lemma \ref{l:bounded-sp}, they have continuous super-potentials and have no mass on proper analytic subsets of $X$. 
Observe that if $S$ is the limit of a sequence $d^{-n_j}(f^{n_j})^*(\phi)$ and $S'$ is a limit of  $d^{-n_j+1}(f^{n_j-1})^*(\phi)$, then $S'$ is in $\Fc$ and $d^{-1}f^*(S')=S$. 
Similarly, if $S''$ is a limit of  $d^{-n_j-1}(f^{n_j+1})^*(\phi)$, then $S''$ is in $\Fc$ and $d^{-1}f^*(S)=S''$. We conclude that $d^{-1}f^*$ maps $\Fc$ to $\Fc$ and this is a surjective map.

Let $\Gc$ denote the projection of $\Fc$ in $H^{p,p}(X,\C)$. This is a compact set in $H^{p,p}(X,\C)$. We deduce from the last discussion that $d^{-1}f^*$ acts on $\Gc$ and is a surjective map. We deduce from Lemma \ref{l:Green+class} that $\Gc$ is a subset of the line $L^+$. So using the above notation, we can write $\{S\}=\lambda_S c^+$ for some real number $\lambda_S$. Using Lemma \ref{l:c+c-} and Proposition \ref{p:Green}, we have
\begin{eqnarray*}
\lambda_S  &=& \lambda_S c^+\smallsmile \{\omega^{k-p}\}=\{S\}\smallsmile \{\omega^{k-p}\}=\lim_{j\to\infty}  \big\langle d^{-n_j}(f^{n_j})^*(\phi),\omega^{k-p}\big\rangle \\
&=& \lim_{j\to\infty}  \big\langle \phi,d^{-n_j}(f^{n_j})_*(\omega^{k-p})\big\rangle = \langle \phi,T^-\rangle =\lambda_\phi.
\end{eqnarray*}
So all currents in $\Fc$ are in the same cohomology class $\lambda_\phi c^+$.

For $S$ in $\Fc$ denote by $\mV_S$ the canonical super-potential of $S-\lambda_\phi T^+$. Recall that $S\leq T^+$ and that the super-potential $\mV_S$ is continuous. Moreover, by Remark \ref{r:main_1}, there is a constant $A>0$ independent of $S$ such that $|\mV_S(R)|\leq A\|R\|_*$ for every $R$ in $\mD_{k-p+1}^0(X)$. 
Recall that $d^{-1}f^*$ is surjective from $\Fc$ to $\Fc$. So for each positive integer $m$, there is a current $S_m$ in $\Fc$ such that $S=d^{-m}(f^m)^*(S_m)$. We have $S-\lambda_\phi T^+=  d^{-m}(f^m)^*(S_m-\lambda_\phi T^+)$ and therefore 
$\mV_S(R)=\mV_{S_m}(d^{-m}(f^m)_*(R))$, according to Proposition \ref{p:pull-current}. 
It follows that 
$|\mV_S(R)|\leq Ad^{-m}\|(f^m)_*(R)\|_*$. Recall that the mass of a positive closed current depends only on its cohomology class and the spectral radius of $f_*$ on $H^{k-p+1,k-p+1}(X,\R)$ is strictly smaller than $d$. We deduce that the RHS of the last inequality tends to 0 as $m$ goes to infinity. Thus, $\mV_S=0$ and hence $S=\lambda_\phi T^+$. This ends the proof of the proposition.
\endproof

\begin{proposition} \label{p:graph}
The sequence of positive closed $(k,k)$-currents $d^{-n}[\Gamma_n]$ converges weakly to the current $T^+\otimes T^-$ in $X\times X$, as $n$ goes to infinity.
\end{proposition}
\proof
Consider a smooth $(r,s)$-form $\phi$ and  a smooth $(k-r,k-s)$-form $\phi'$ on $X$. 
By Weierstrass approximation theorem, the vector space spanned by forms of type $\Phi:=\phi(y)\wedge \phi'(x)$ is dense in the space of $(k,k)$-forms on $X\times X$. Therefore, it is enough to check that $\langle d^{-n} [\Gamma_n], \Phi\rangle$ converges to $\langle T^+\otimes T^-, \Phi\rangle$.
Note that
$\langle d^{-n} [\Gamma_n], \Phi\rangle = \langle d^{-n} (f^n)^*(\phi),\phi' \big\rangle$. By this identity and
Lemma \ref{l:conv-rs}, we obtain that $\langle d^{-n} [\Gamma_n], \Phi\rangle$ tends to 0 if $(r,s)\not=(p,p)$.
Observe also that
$\langle T^+\otimes T^-, \Phi\rangle = \langle T^-,\phi\rangle\langle T^+,\phi'\rangle$
and the last product vanishes unless $(r,s)=(p,p)$. So the proposition is clear when $(r,s)\not=(p,p)$.

Assume now that $r=s=p$. By Proposition \ref{p:Green-eq}, 
$\langle d^{-n}[\Gamma_n],\Phi\rangle$, which is equal to $\big\langle d^{-n} (f^n)^*(\phi),\phi' \big\rangle$,
converges to $\langle T^-,\phi\rangle\langle T^+,\phi'\rangle$. We have seen that the last product is equal to $\langle T^+\otimes T^-,\Phi\rangle$. The proposition follows.
\endproof

In order to prove Theorem \ref{t:main_4}, we need the following lemma which was obtained in \cite[Prop. 2.2]{DinhNguyenTruong16}.

 \begin{lemma} \label{l:DNT}
 Let $V_n$ be a sequence of effective $k$-chains in $X\times X$ for $n\in\N$.  Suppose there is a sequence of positive numbers 
 $d_n$ converging to infinity such that $d_n^{-1}[V_n]$
 converges weakly to a positive closed $(k,k)$-current $R$ on $X\times X$.
 Assume also that the tangential  h-dimension of $R$ with respect to $\Delta$ is minimal, i.e. equal to $0$. Define 
 $c:= \{R\}\smallsmile \{\Delta\}$. Then the number 
$\delta_n$ of isolated points in the intersection of $V_n$ with $\Delta$, counted with multiplicity, satisfies $\delta_n\leq cd_n+o(d_n)$ as $n$ tends to infinity. 
 \end{lemma}

\noindent
{\bf End of the proof of Theorem \ref{t:main_4}.}
We apply Proposition \ref{p:graph} and the last lemma for $\Gamma_n, T^+\otimes T^-, d^n, P_n$ instead of $V_n, R, d_n,\delta_n$.
Since $T^+$ has a continuous super-potential, the wedge-product $T^+\wedge T^-$ is well-defined. Therefore, by Theorem \ref{t:main_2}, the tangential h-dimension of $T^+\otimes T^-$ along $\Delta$ is 0. We also have 
$$\{T^+\otimes T^-\}\smallsmile \{\Delta\}=\{T^+\}\smallsmile\{T^-\}=c^+\smallsmile c^-.$$ 
The last cup-product is 1 by Lemma \ref{l:c+c-}. Applying Lemma \ref{l:DNT} gives the result. 
\hfill $\square$

%%%%%%%%%%%%%%%%%%%%%%%%%%%%%%%%%%%%%%%%%%%%%%%%%%%%%%%%%%%%%%%%%


\small

\begin{thebibliography}{99}

\bibitem{BedfordTaylor}
Bedford E., Taylor B.A., A new capacity for plurisubharmonic
functions. {\it Acta Math.} {\bf  149}  (1982), no. 1-2, 1-40. 

\bibitem{Blanchard} 
Blanchard A., Sur les vari\'et\'es analytiques complexes. 
{\it Ann. Sci. \'Ecole Norm. Sup.} (3) {\bf 73} (1956), 157-€"202. 

\bibitem{BGS}
Bost J.-B., Gillet H., Soul\'e C., Heights of projective varieties and positive Green forms. {\it J. Amer. Math. Soc.} {\bf 7} (1994), no. 4, 903-1027.

\bibitem{Demailly92} 
Demailly J.-P., Courants positifs et th\'eorie de l'intersection.  {\it Gaz. Math. No.}  {\bf 53} (1992), 131-159. 

\bibitem{Demailly12} 
Demailly J.-P.,
 {\em Complex Analytic and Differential Geometry},   2012. Available online at 
{\tt www-fourier.ujf-grenoble.fr/$\sim$demailly/books.html}


\bibitem{Demailly_survey}
 Demailly J.-P.,  Analytic methods in algebraic geometry. {\it Surveys of Modern Mathematics} {\bf 1} (2012), International Press, Somerville, MA; Higher Education Press,
              Beijing.


\bibitem{Diller_Favre}
Diller  J., Favre  C.,  Dynamics of bimeromorphic maps of surfaces. {\it Amer. J. Math.} {\bf 123} (2001), no. 6, 1135-1169.     

 \bibitem{Dinh05}
Dinh  T.-C.,   Suites d'applications m\'eromorphes multivalu\'ees et courants laminaires. {\it J. Geom. Anal.}  {\bf 15} (2005),
207-227.

\bibitem{DinhNguyen06}
Dinh T.-C.,  Nguyen V.-A.,  The mixed Hodge-Riemann bilinear relations for compact 
K{\"a}hler manifolds. {\it Geom. Funct. Anal.}, Vol. {\bf 16} (2006), 838-849. 


 \bibitem{DinhNguyenTruong15} 
 Dinh T.-C., Nguyen V.-A., Truong T.-T.,
  Equidistribution for meromorphic maps with dominant topological degree.
 {\it  Indiana Univ. Math. J.}  {\bf 64} (2015), no. 6, 1805-1828.
   


\bibitem{DinhNguyenTruong16}
Dinh T.-C,  Nguyen V.-A., Truong T.-T., Growth of the number of periodic points for meromorphic maps. {\it Bull. London Math. Soc.}, to appear. 
 {\tt  DOI: 10.1112/blms.12082}


\bibitem{DinhSibony04} 
Dinh T.-C., Sibony N.,
Regularization of currents and entropy. {\it Ann. Sci. \'Ecole Norm. Sup. (4)}  {\bf  37} (2004), no. 6, 959-971.



\bibitem{DinhSibony07}
  Dinh T.-C., Sibony N., Pull-back of currents by holomorphic maps.
{\it  Manuscripta Math.} {\bf 123} (2007), 357-371.

\bibitem{DinhSibony09}
 Dinh T.-C., Sibony N.,
 Super-potentials of positive closed currents, intersection theory and dynamics. {\it Acta Math.}  {\bf 203} (2009), no. 1, 1-82.


\bibitem{DinhSibony10}
Dinh T.-C., Sibony N., Super-potentials for currents on compact K\"ahler manifolds and dynamics of automorphisms.
{\it J. Algebraic Geom.} {\bf  19} (2010), no. 3, 473-529. 



\bibitem{DinhSibony12}
Dinh T.-C., Sibony N., 
Density of positive closed currents, a theory of non-generic intersections. {\it J. Algebraic Geom.}, to appear.
{\tt arXiv:1203.5810}  


\bibitem{DinhSibony15}  
Dinh T.-C., Sibony N.,  Unique ergodicity for foliations in $\P^2$ with an invariant curve. {\it Invent. Math.}, to appear.
{\tt DOI 10.1007/s00222-017-0744-2} 

\bibitem{DinhSibony16b}
Dinh T.-C., Sibony N.,
 Equidistribution problems of complex dynamics in higher dimension. {\it Int. J. Math.} {\bf 28} (2017), no. 7, 1750057, 31 pp. 
 
 
 \bibitem{Favre} 
 Favre C.,  Points p\'eriodiques d'applications birationnelles de $\P^2.$ {\it Ann. Inst. Fourier (Grenoble)} {\bf  48} (1998), no. 4, 999-1023.
 
 
\bibitem{FornaessSibony}
Forn\ae ss J.-E., Sibony N., Oka's inequality for currents and
applications. \textit{Math. Ann.} \textbf{301} (1995), 399-419.
 
 
\bibitem{Gromov}
Gromov M.,  Convex Sets and K\"ahler Manifolds. World Sci. Publishing, Teaneck, NJ, 1990.


\bibitem{IwasakiUehara}  
 Iwasaki K., Uehara T., Periodic points for area-preserving birational maps of surfaces. 
{\it Math. Z.}  {\bf 266} (2010), no. 2, 289-318.  


\bibitem{Oguiso}
Oguiso K., 
A few explicit examples of complex dynamics of inertia groups on surfaces - a question of Professor Igor Dolgachev. {\it Preprint} (2017). {\tt  arXiv:1704.03142}


 \bibitem{Saito}   
 Saito S.,  General fixed point formula for an algebraic surface and the theory of Swan representations for two-dimensional local rings.  {\it Amer. J. Math.} {\bf 109} (1987), no. 6, 1009-1042.   


\bibitem{Skoda}
Skoda H.,  Prolongement des courants positifs, ferm{\'e}s de
masse finie. \textit{Invent. Math.} \textbf{66} (1982), 361-376.


 


\bibitem{Vu}  
Vu D.-V., Intersection of positive closed currents of higher bi-degree. 
{\it Michigan Math. J.}  {\bf  65} (2016), no. 4, 863-872.
 
  
 \bibitem{Xie} 
 Xie J., Periodic points of birational transformations on projective surfaces.
{\it Duke Math. J.}  {\bf 164} (2015), no. 5, 903-932. 


 
\end{thebibliography}
\end{document}